\documentclass[ssy,preprint]{imsart}
\usepackage{natbib}
 \bibpunct[, ]{(}{)}{,}{a}{}{,}%
 \def\bibsep{\smallskipamount}%
\arxiv{arXiv:0000.0000}

\usepackage{graphicx,color,subfigure}
\usepackage[displaymath,mathlines,pagewise]{lineno}
\usepackage{algorithm,algorithmicx,algpseudocode}
\usepackage{amsmath,amstext,amsfonts,amssymb,mathrsfs,euscript}
\usepackage{personal}

\begin{document}
\begin{frontmatter}
\title{Transitory Queueing Networks}
\runtitle{Transitory Networks}
\begin{aug}
\author{\fnms{Harsha} \snm{Honnappa}\ead[label=e1]{honnappa@purdue.edu}}
\address{School of Industrial Engineering,\\ Purdue University,\\ West
  Lafayette IN 47906.\\\printead{e1}\\},
\author{\fnms{Rahul} \snm{Jain}\ead[label=e2]{rahul.jain@usc.edu}}
\address{EE, CS \& ISE Departments,\\ University of Southern
California,\\ Los Angeles, CA 90089. \printead{e2}}
\affiliation{Purdue University\thanksmark{m1} and University of
  Southern California\thanksmark{m2}}
\runauthor{Honnappa and Jain}
\end{aug}
% \date{Received: \today}

\runauthor{Honnappa and Jain}

\runtitle{Transitory Networks}

% \TITLE{Bundling Information Goods of Decreasing Value}
%\TITLE{Transitory Queueing Networks}

% \ARTICLEAUTHORS{%
% \AUTHOR{Harsha Honnappa}
% \AFF{School of Industrial Engineering, Purdue University, IN 47906 \EMAIL{honnappa@purdue.edu}}
% \AUTHOR{Rahul Jain}
% \AFF{Department of Electrical Engineering, University of Southern
%   California, CA 90089, \EMAIL{rahul.jain@usc.edu}}
 % end of the block

\begin{abstract}
Queueing networks are notoriously difficult to analyze \textit{sans}
both Markovian and stationarity assumptions. Much of the theoretical
contribution towards performance analysis of time-inhomogeneous
single class queueing networks has focused on Markovian networks, with
the recent exception of work in \cite{LiWh2011} and \cite{MaRa10}. In this paper,
we introduce \textit{transitory queueing networks} as a model of
inhomogeneous queueing networks, where a large, but finite, number of
jobs arrive at queues in the network over a fixed time horizon. The
queues offer FIFO service, and we assume that the service rate can be
time-varying. The non-Markovian dynamics of this model complicate the
analysis of network performance metrics,
necessitating approximations. In this paper we develop fluid and
diffusion approximations to the number-in-system performance metric by
scaling up the number of external arrivals to each queue, following
\cite{HoJaWa2012}. We also discuss the implications for bottleneck
detection in tandem queueing networks.
\end{abstract}
% \KEYWORDS{transient analysis; fluid approximations; diffusion
%   approximations; bottlenecks}
% {\em OR/MS subject classification:} Queues:Networks, Approximations,
% Probability:Diffusion, Production/Scheduling:Stochastic \\
% %\subclass{MSC 91A10 \and MSC 91A80 \and MSC 91B26.}
% \HISTORY{This paper was first submitted July 18, 2017.}
\end{frontmatter}
\section{Introduction} \label{intro}
Single class queueing networks (henceforth `queueing networks') have been studied extensively in the
literature, with much effort focused on understanding the steady-state
joint distribution of the state of the network (typically defined as
the number of jobs in each queue). In this paper, we consider a
variation of the generalized Jackson queueing network model
(\cite[Chapter 7]{ChYa01}) where a finite, but large, number of jobs
arrive at some of the nodes in the network from an extraneous source. We characterize fluid and
diffusion approximations to the queue length state process, as the
population size scales to infinity. %The motivations for this model
%include transportation networks, manufacturing and service networks.  % To the best of our knowledge, this
% is the first time transitory networks have been studied in the
% literature, despite its wide applicability.

To motivate this model, consider a manufacturing facility that produces engines. Each part of the engine is
produced and assembled in a separate machine that requires some human
supervision, with final assembly occuring at the end of the job-shop. Typically, there is a finite, but large number of jobs that need to
be completed in a shift spanning a few hours. Furthermore, jobs cannot be carried over to the
next shift. It is typically the case that the shift horizon is not
long enough for the system to reach a steady state. Furthermore, 
the bottleneck machine can change over the shift horizon, as a
consequence of the variation in the arrival of jobs to the job-shop,
or due to stochastic variation in the assembly time in each
machine. This is also known as the `shifting bottleneck'
phenomenon in production engineering; see \cite{RoNaTa2002,LaBu1994,LaBu1995}. In
the latter two papers, Jackson networks are used to analyze the shifting
bottleneck network. More generally, there is less reason to believe
that the arrival and service processes are stationary and ergodic, and
that they are Poisson processes. We introduce \textit{transitory
  queueing networks} as a broadly applicable model of such systems. 

Transitory queueing networks consist of a number of infinite buffer,
FIFO, single server queues (a.k.a. `nodes') interconnected by customer
routes. We assume that the routing matrix satisfies a so-called 
\emph{Harrison-Reiman (H-R)} condition that the matrix has a spectral
radius of less than one. This implies that, on
completion of service at a particular node, a customer is routed to
another node or exits the network altogether. 
% Here, we are interested in a variation of such networks, called 'transitory' generalized Jackson
% networks. These are networks of queues where the total
% number of jobs that arrive for service is finite, or the time horizon
% of interest is finite. Such networks are models of
% service and communication networks that operate over a finite time
% horizon, or the system operator is interested in short-time behavior
% of such networks. Further, traffic and service in most queueing
% networks exhibit time inhomogeneties. For example, in datacenters and hospitals demand exhibits very
% obvious time-of-day and day-of-week effects. In these circumstances,
% standard analyses of the performance of generalized Jackson networks
% will not suffice, and one must resort to a purely transient, or 'transitory' analysis of the network.
Jobs enter a given node at random time epochs modeled as the ordered statistics of independent and
identically distributed (i.i.d.) random variables. The arrival epoch
of the jobs at different nodes can be correlated (for instance,
in the manufacturing context, jobs can be submitted to multiple nodes simultaneously). We assume that the service processes
at different nodes are independent with time-inhomogeneous service
rates, and  modeled as a time change of a unit rate
renewal counting process, generalizing the construction of a
time-inhomogeneous Poisson process. The transient analysis of generalized Jackson networks is non-trivial,
as noted before. The conventional heavy-traffic diffusion
approximation that relied on long-run average rates has been used to approximate the evolution of the state
process. However, these rates do not exist in transitory queueing
networks. Here, we develop a `population
acceleration' approximation, by
increasing the number of jobs arriving at the network in the interval of
interest to infinity, and suitably scaling (or `accelerating') the
service process in each queue by the population size. 

Our results complement
the existing literature on the analysis of single-class queueing
networks by establishing the following results: 

\begin{enumerate}
\item[(i)] we develop a large population approximation framework for
  studying single class queueing networks in a transitory setting,
  complementing and extending Markovian network analyses to 
  non-Markovian queueing networks, 
  \item[(ii)] Our diffusion approximations
  use the recently developed directional derivative oblique reflection
  map in \cite{MaRa10} to establish a diffusion scale approximation;
  this is substantially different from the conventional heavy-traffic
  approximations used to study single-class queueing networks, and
  \item[(iii)] we study the evolution of the bottleneck process over the time
  horizon, identifying the bottleneck station as time progresses. This
  analysis extends the standard bottleneck analyses, where bottlenecks
  are identified in terms of the long-term average arrival and service
  rates.
\end{enumerate}

\subsection{Analytical Results}
We consider a
sequence of queueing networks wherein
$n$ jobs arrive at each node that receives external traffic in the
$n$th network. We first establish a
functional strong law of large numbers (FSLLN) to the arrival
process, as the population size
scales to infinity, by generalizing the Glivenko-Cantelli Theorem to
multiple dimensions.  We also assume the service processes
satisfy a FSLLN and functional central limit theorem (FCLT) in the population acceleration scale. The queue length fluid limit is shown to be equal to the oblique
reflection of the difference of the fluid arrival and service
processes (or the `fluid netput' process). On the other hand, for the FCLT we introduce the
notion of a multidimensional Brownian bridge process, as a
generalization of a the one-dimensional Brownian bridge process, and show that
the diffusion scaled arrival process converges to a multidimensional
Brownian bridge in
the large population limit. The diffusion limit turns out to be
complicated, and it is shown to be a reflection of a multidimensional
Gaussian bridge process - however, the reflection is through a
directional derivative of the oblique reflection of the netput in the
direction of the diffusion limit of the fluid netput process. This is a
highly non-standard result. Indeed, it is only in the recent past that
\cite{MaRa10} have investigated the existence
of a directional derivative to the oblique reflection map.

Leveraging the results of \cite{MaRa10}, we can only establish a pointwise
diffusion limit for an arbitrary transitory
queueing network. This is due to the fact that the
directional derivative limit can have sample paths with
discontinuities that are both right- \textit{and} left-discontinuous. Thus,
establishing convergence in a sample path space under a suitably weak
topology such as the (`strong' or `weak') $M_1$, for instance, is not 
straightforward. Instead, we focus on the case of tandem queueing networks, with
uniform and unimodal arrival time distribution functions. In this
case, we show that the discontinuities in the limit are either right
or left continuous, and hence we can establish $M_1$ convergence. Using these approximations, we next address the question of
bottleneck prediction. We generalize the standard definition of a
bottleneck in a single class network, defined as the queues whose
fluid arrival rate exceeds the fluid service capacity to the transitory
setting.%  We focus primarily on the departure time of a job from a
% tandem queueing network. 
% There is an extensive literature focused on the
% analysis of single class queueing networks.

\subsection{Implications for Bottleneck Detection}
Bottleneck detection and prediction is likely the most important
question that a system operator faces. Heavy traffic theory has been
immensely successful at characterizing steady state bottlenecks in 
general queueing networks, under minimal network data
assumptions. However, there are many circumstances, ranging from
manufacturing, to healthcare, transportation and computing, where transient
bottleneck detection and analysis is critical. In the purely
transient, or `transitory' setting, it is common for the bottleneck node to change over
the shift horizon. As a consequence the plant manager moves
workers around trying to ease bottlenecks, increasing costs and
increasing the likelihood of job overages. Another example is in 
healthcare where patient diagnosis relies on a
number of tests that must be done with different
machines. Furthermore, in many time critical settings, the horizon
within which tests must be conducted is fixed. An important question
in these situations is
whether transitory bottlenecks can be accurately predicted, given network data.

Note that
the standard definition of a bottleneck is a `capacity level' one,
defined in terms of long-term averages. This definition, of course, is
not satisfactory in the transitory setting where steady states might
not be reached. Second, standard
heavy-traffic analysis completely ignores non-stationarities in the
network data, which is generally prevalent in the examples presented above, making it
an inappropriate analytical tool to use. In particular, it is possible
that the number of bottleneck nodes in the network changes over time,
thus necessitating a more nuanced definition and analysis of
bottlenecks. We introduce a natural definition of a `transitory bottleneck' in
terms of the diffusion approximations developed in this paper, and then use
these approximations to analyze the evolution of the number of
bottlenecks in tandem transitory networks. This is particularly
relevant to job-shop analysis of production systems.

\subsection{Related Literature}
There has been significant interest in the analysis of single class queueing
networks. Under the assumption of Poisson arrival and service processes, \cite{Ja1957} showed that the steady state distribution of the state
of the network (the number of jobs waiting in each node) is equal to
the product of the distribution of the state of each node in the
network. This desirable property implies that, in steady state, the
network exhibits a nice independence property. This property does not
extend to networks with general arrival and service processes; these
are also known as generalized Jackson networks. 

Reiman first established the heavy-traffic diffusion approximation to open
generalized Jackson networks in \cite{Re1984}. In particular, the
diffusion approximation is shown to be a multi-dimensional reflected
Brownian motion in the non-negative orthant, reflected through the
oblique reflection mapping. Such reflection maps have come to be
called as Harrison-Reiman maps following the seminal work in
\cite{HaRe1981}. \cite{ChMa1991a,ChMa1991}
characterize a homogeneous fluid network, and establish
fluid and diffusion approximations. The analysis of non-stationary and
time inhomogeneous queueing systems is non-trivial in general. For
single server queues, see \cite{Ke82,Ma85,MaMa95} among others. In
\cite{HoJaWa2012,HoJaWa2013b} we develop fluid and diffusion
approximations to the $\D_{(i)}/GI/1$ transitory single server
queue in the population acceleration scaling regime. Recent work in \cite{BeHoJvL2015,BeHoJvL2016} considers the
$\D_{(I)}/GI/1$ queue under a `uniform' acceleration regime, where the
initial work in the network is assumed to satisfy a heavy-traffic like
condition. For networks of queues,
\cite{MaMaRe1998} develops strong approximations to queueing networks
with nonhomogeneous Poisson arrival and service processes. In
\cite{DuMaWh2001}, the authors study the offered load process in a
bandwidth sharing network, with nonstationary traffic and general
bandwidth requirements. The closest work in the literature to
the current paper is \cite{LiWh2011} who
study a network of non-Markovian fluid queues with time-varying
traffic and customer abandonments. To be precise, they consider a
$(G_t/M_t/s_t + GI_t)^m/M_t$ network with $m$ nodes, time-varying
arrivals, staffing and abandonments, and inhomogeneous Poisson service
and routing, and characterize the performance of the network as a
direct extension of the single-server queue case. However, their
scaling regime and the limit processes are completely different from
the our results.

The rest of the paper is organized as follows. We collect relevant
notation in Section~\ref{sec:notation}. We start our analysis with a
description of the transitory generalized Jackson network model in
Section \ref{sec:primitives}, and we develop fluid and diffusion
approximations to the network primitives. In Section
\ref{sec:fluid-networks}, we develop functional strong law of large
numbers approximations to the queueing equations, and identify the
fluid model corresponding to the transitory network. We identify the
diffusion network model in Section
\ref{sec:diffusion-networks}, and establish a weak convergence result
for a tandem network with unimodal arrival time distribution. We end
with conclusions and future research directions in Section \ref{sec:conclusions}

\section{Notation}\label{sec:notation}
Following standard notation, $\sC^K$ represents the space of continuous
$\bbR^K$-valued functions, and $\sD^K$ the space of functions that are
right continuous with left limits (RCLL) and are $\bbR^K$-valued. The space
$\sD^K_{l,r}$ consists of $\bbR^K$-valued functions that are either right- or
left-continuous at each point in time, while $\sD^K_{\lim}$ is the space
of $\bbR^K$-valued functions that have right and left limits at all
points in time. $\sD^K_{usc}$ is the space of
RCLL functions that are upper semi-continuous as well. We represent
the space of $L\times L$ matrices by $\sM_L$. The space
and mode of convergence of a sequence of stochastic elements is
represented by $(X,Y)$, where $X$ is the space in which the stochastic
elements take values and $Y$ the mode of convergence. In this paper
our results will be proved under the uniform mode of convergence and
occasionally in the ``strong'' $M_1$ ($SM_1$) topology (see \cite[Chapter
11]{Wh01}). Weak convergence of measures will be represented by
$\Rightarrow$. Finally, $\text{diag}(x_1, \ldots, x_K)$ represents a $K
\times K$ diagonal matrix with entries $x_1,\ldots,x_K$.
\section{Transitory Queueing Network} \label{sec:primitives}
% We consider a single class queueing network with $K$ single server FIFO
% nodes. Each node starts service at some fixed time, which could be
% different from the other nodes. The servers are
% non-preemptive and non-idling. We assume that a finite, but large,
% number ($n$) of jobs apply for service at each node, and that every job
% is served independently of the others. The network is assumed to offer
% Markovian routing between the nodes. Thus, the routing can be
% represented by a sub-stochastic routing matrix $\mathbf P$. Furthermore, we assume that the network is open implying that all
% arriving users eventually depart the network. In this section we
% present (and prove, where necessary) functional strong law of large numbers and functional central limit theorem results for the
% network data; that is, the arrival process $\mathbf A$, the service
% process $\mathbf S$ and the routing process $\mathbf R$, in the limit
% of a large number of arrivals $n$ by rescaling the service process
% appropriately by the population size. We call this the `population
% acceleration' approximation regime, analogous to the `uniform
% acceleration' regime used in \cite{MaMaRe1998}.

Let $(\Omega, \mathcal{F}, \mathbb P)$ be an appropriate probability space on
which we define the requisite random elements. Let $\sK :=
\{1,\ldots,K\}$ be the set of nodes in the network, and $\sE
\subset \sK$ the set of nodes where exogenous traffic enters the
network. Each node in $\sE$ receives $n$ jobs that arrive exogenously
to the node. We assume a very general model of the traffic: let $\mathbf
T_m := (T_{e_1,m},\ldots,T_{e_J,m})\,\,, m \leq n,$ represent the tuple of
arrival epoch random variables where $T_{e_j,m}$ is the arrival epoch of
the $m$th job to node $e_j \in \sE$ (here $J :=
|\sE|$). By assumption $T_{e_j,m} \in [0,T]$ for all $e_j \in \sE$ and $1
\leq m \leq n$. We also assume that $\{\mathbf T_m ; m =
1,\ldots,n\}$ forms a sequence of independent random vectors. Let
$F_{e_j}$ be the distribution function of the arrival epochs to node $e_j
\in \sE$; that is $\bbE[\mathbf 1_{\{T_{e_j,m} \leq t\}}] = F_{e_j}(t)$ with
support $[0,T]$. Users sample a time epoch to arrive at the node and enter the
queue in order of the sampled arrival epochs; thus the arrival process
to each node is a function of the ordered statistics of the arrival
epoch random variables. In many situations, it is plausible
that there is correlation between the arrival processes to the
nodes in $\sE$. To model such phenomena, we assume that the joint
distribution of the arrival epochs are
fully specified. To be precise, we assume that $\bbP(T_{e_1,m} \leq t,
\ldots, T_{e_J,m} \leq t)$ for all $m \in \{1,\ldots,n\}$ is well
defined. Let $\mathbf a_m(t) := (\mathbf 1_{\{T_{e_1,m} \leq t\}},
\ldots, \mathbf 1_{\{T_{e_J,m} \leq t\}}) \in \sD^{J}[0,\infty)$ and
\begin{equation} \label{eq:arrival process}
A_{n,e_j} := \sum_{m=1}^n \mathbf 1_{\{T_{e_j,m} \leq t\}} ~\text{for}~ 1 \leq j
\leq J,
\end{equation}
then $\mathbf A_n(t) := \sum_{m=1}^n \mathbf a_m(t) =
(A_{n,e_1}(t), \ldots, A_{n,e_J}(t)) \in \sD^{J}[0,\infty)$ is the vector of
cumulative arrival processes to the nodes in $\sE$. Then,
\[
\bbE[\mathbf A_n(t)] = n \mathbf F(t) := n (F_{e_1}(t), \ldots,
F_{e_J}(t))
\]
and $\bbE[\mathbf A_n(t) \mathbf A_n(t)^T] = [n F_{e_i,e_j}(t) + n(n-1) F_{e_i}(t) F_{e_j}(t)]$,
where $F_{e_i,e_j}(t) := \bbP(T_{e_i,m} \leq t, T_{e_j,m}\leq t)$. This `multivariate empirical
process' representation for the traffic affords a very
natural model of correlated traffic in networks, and stands in
contrast with generalized Jackson networks where external traffic to
each node in $\sE$ is independent.

Recall from
Donsker's Theorem (for empirical sums) that $n^{-1/2} (A_{n,e_i} - n F_{e_i})
\Rightarrow W_{e_i}^0 \circ F_{e_i}$, where $W_{e_i}^0$ is a standard Brownian bridge
process. The Brownian bridge process $W^0$ is also well defined as a `tied-down'
Brownian motion process equal in distribution to $(W(t) - t
W(1),~t \in [0,1])$, for all $t \in [0,1]$, where $W$ is a standard
Brownian motion process.%  Assume that the covariance matrix of
% $\mathbf T_m$ is $R$.

\begin{definition} \label{def:bridge}
  Let $\mathbf W = (W_{1},\ldots, W_{J})$ be a $J$-dimensional standard Brownian motion
  process with identity covariance matrix. If $R$ is a $J\times J$ positive-definite matrix with 
  lower-triangular Cholesky factor $L$, then $\mathbf W_R = L \mathbf W$ is a
  $J$-dimensional Brownian motion with covariance matrix $R$. By directly
  extending the definition of a one-dimensional Brownian bridge
  process,
  \[
  \left( \mathbf W^0(t) = \mathbf W_R(t) - t \mathbf W_R(1),~t \in
    [0,1]\right)
  \]
  is a $J$-dimensional Brownian bridge process with covariance matrix $R$.
\end{definition}

It is straightforward to see that $\mathbb E[\mathbf W^0(t)] = 0$ for
all $t \in [0,1]$ and $\mathbb E[\mathbf W^0(t) \mathbf W^0(s)^T] =
t(1-s) R \equiv [t(1-s) r_{i,j}]$ when $t \leq s$. More generally, we
define a Brownian bridge process with time-dependent covariance matrix
as follows. Recall that a stochastic process is defined as a Gaussian
process provided its finite dimensional distributions are jointly Gaussian.

\begin{definition}~\label{def:bb-multi}
  Let $R : [0,T] \times [0,T]\to \sM_J$ be a right continuous, symmetric function such that $R(t,s)$ is
  positive-definite for each $t,s \in [0,T]$, with the restriction
  that $R(0,0) = 0$ and $r_{i,j}(T,T) = r_i(T)r_j(T)$, where $r_{i,j}$ is the $(i,j)$th
  entry of $R$ and $r_i$ is the $i$th
  diagonal element of $R$. The Gaussian process $\mathbf W^0$ is a  
  $J$-dimensional Brownian bridge process if it has mean zero and
  covariance function $\bbE[\mathbf W^0(t)
  \mathbf W^0(s)^T] = R(t,s) - \text{diag}(R(t,t)) \text{diag}(R(s,s))
  \equiv [r_{i,j}(t,s) - r_i(t) r_j(s)]$, where $\text{diag} : \sM_J
  \to \bbR^J$ is a function that maps a matrix to a vector of the
  diagonal elements.
\end{definition}

Note that this definition is a natural generalization of the bridge
process in Definition~\ref{def:bridge}. The terminal condition
$r_{i,j}(T,T) = r_i(T)r_j(T)$ ensures that $\mathbf W^0(T) = 0$
a.s. While we do not argue the existence of this object rigorously, the right
continuity of the covariance matrix and the assumed Gaussianity of the
marginals imply that it can be inferred from the Daniell-Kolmogorov
theorem; see~\cite{RoWi1988}. Now, observe that the covariance
function of the pre-limit traffic process satisfies $R_n(t,s) =
\bbE[(\mathbf A_n(t) - \bbE[\mathbf A_n(t)])(\mathbf
A_n(s) - \bbE[\mathbf A_n(s)])^T] = n [F_{e_i,e_j}(t,s) - F_{e_i}(t)
F_{e_j}(s)]\in \sC^{J \times J}$, where $F_{e_i,e_j}(t,s) :=
\bbP(T_{e_i,m} \leq t, T_{e_j,m} \leq s)$ and $t \leq s$. Now,
following Definition~\ref{def:bb-multi} let
$\mathbf W^0 \circ \mathbf F$ represent a multidimensional Brownian
bridge with covariance function
\begin{align}
  \label{eq:bb-cov}
  R(t,s) = [F_{e_i,e_j}(t,s) - F_{e_i}(t) F_{e_j}(s)],
\end{align}
when $t \leq s$. Note that we are ``overloading'' the composition
operator $\circ$ in this notation, but the usage should be clear from
the context. Theorem~\ref{thm:arrival-limits} below establishes
multivariate generalizations of the classical Glivenko-Cantelli and Donsker's
theorems.

\begin{theorem} \label{thm:arrival-limits}
% Consider the triangular array of i.i.d. random vectors $\{\mathbf
% T_m,~m \leq n \}$ $n \geq 1$, and let $\mathbf a_m(t) := (\mathbf 1_{\{T_{e_1,m} \leq t\}}, \ldots, \mathbf 1_{\{T_{e_J,m} \leq t\}})$, for
% $t \in [0,\infty)$. Then,
We have,

\noindent (i) \(
n^{-1} \mathbf A_n \to \mathbf F  
\)
in $(\sC^J, U)$ a.s. as $n \to \infty$, and

\noindent (ii) 
\(
\hat{\mathbf A}_n := n^{-1/2}\left( \sum_{m=1}^n \mathbf a_m - n \mathbf F \right) \Rightarrow \mathbf W^0 \circ \mathbf F
\)
in $(\sC^J, U)$ as $n \to \infty$, where $\mathbf W^0 \circ \mathbf F\in
\sC^J[0,\infty)$ is the $J$-dimensional Brownian bridge process with
covariance matrix defined in~\eqref{eq:bb-cov}.
\end{theorem}
The proof of the theorem is available in the appendix. We refer to the
$k$th component process by $W_k\circ F_k$.

Next, we consider a sequence of service processes indexed by
the population size $n \geq 1$, $S_{n,k} :
\Omega \times [0,\infty) \to \bbN$ for $k \in \sK$. We assume that for each $k \in \sK$ the function
$\m_{n,k} : [0,\infty) \to [0,\infty)$ is
Lebesgue-integrable and that $M_{n,k}(t) := \int_0^t \m_{n,k}(s) ds$ satisfies $ M_{n,k} \to M_k
~\text{in}~(\sC,U)~\text{as}~ n \to \infty$, where $M_k : [0,\infty) \to [0,\infty)$ is non-decreasing and continuous. We also assume that
\begin{equation}
  \label{eq:1}
  \mathbf M_n := (M_{1,n}, \ldots, M_{n,k}) \to \mathbf M
  := (M_1, \ldots, M_K) ~\text{in}~(\sC^K,U)~\text{as}~ n \to \infty,
\end{equation}
where $K = |\sK|$.
% For instance, a natural scaling would be to
% assume that the
% $\{\nu_m^{n,k},~m \geq 1\} = \{\nu_m^k/n,~m\geq 1\}$, where
% $\{\nu_m^k,~m\geq 1\}$ form an
% i.i.d. sequence of random variables.  % such that $E S_{n,k}(t) = M^n_k(t)$ for any $t \geq 0$. 
% It is useful to think of $\mu_{n,k}(t)$ as the `service rate function' associated with a
% system with $n$ arriving jobs.
Let $\mathbf S_n := (S_{1,n}, \ldots,
S_{n,k})$ represent the `network' service process, where the component service
processes are independent of each other. We
assume that $\mathbf S_n$ satisfies the following fluid and diffusion
limits. 

\begin{assumption} \label{thm:service-limits}
The service processes $\{\mathbf S_n, ~n\geq 1\}$ satisfies

\noindent (i) 
\(
\left[n^{-1} \mathbf S_n- \mathbf M_n \right] \to 0 
\)
in $(\sC^K,U)$ a.s. as $n \to \infty$, and

\noindent (ii)
\(
\hat {\mathbf S}_n(t) := n^{-1/2} \left( \mathbf S_n - n \mathbf M \right) \Rightarrow \mathbf W \circ \mathbf M
\)
in $(\sC^K,U)$ as $n \to \infty$, where $\mathbf W := (W_1, \ldots,
W_K)$ is a $K-$dimensional Brownian motion process with identity covariance
matrix.
\end{assumption}

Note that the covariance function of the process $\mathbf W \circ
\mathbf M$
is the diagonal matrix with entries $(M_1,\ldots, M_K)$. This service process is analogous to the time-dependent
`general' traffic process $G_t$ proposed in
\cite{LiWh2014}. It's possible to anticipate a proof of this
result when the centered service process $\mathbf S_n - \mathbf M_n$
is a martingale. This would be the case when $\mathbf S_n$ is a
$K$-dimensional stochastic process where the marginal processes are
nonhomogeneous Poisson processes and $M_{n,k} = \bbE[S_{n,k}]$. Here, we
leave the development of a general result to a separate paper
and, instead, assume that such a sequence of service processes exist.

% The service times of each user are assumed to be independent of the
% others'. Let $\nu_i^k : \Omega \rightarrow (0,\infty)$ represent the
% service time of the $i$th arrival to the $k$th node in the
% network. % Let $m_k = 1/\mu_k$ be the mean service time and $c_k^2$ the
% % (finite) squared coefficient of variation of the service times. 
% % Suppose that service starts at time $T_{s,k} \geq 0$ in the $k$th node. Then, t
% The number of service completions by time $t$ at node $k$ is $S_k(t)
% := sup \{m \leq n | \sum_{i=1}^m \nu_i^k \leq t\}$. The tuple $\mathbf
% S(t) := (S_1(t), \ldots, S_K(t)) \in \sD^K[0,\infty)$ represents the
% service processes in the network. By assumption, the components of
% this vector are independent. We assume that for each server $k \in \sK$ there exists a
% non-negative continuous function $\mu_k : [0,\infty) \to [0,\infty)$
% such that mean service process satisfies $\bbE[S_k(t)] = M(t) :=
% \int_0^t \mu_k(s) ds$. We also assume that the covariance
% function $\bbE[(S_k(t)-
% \bbE[S_k(t)])(S_k(t)-\bbE[S_k(t)])^T]$ exists and is finite for each
% $k \in \sK$.

On completion of service at node $i$, a job will join node $j$ with
probability $p_{i,j} \geq 0 \,\, , \,\,i,j \in \{1,\ldots,K\}$, or
exit the network with probability $1-\sum_{j} p_{i,j}$. Thus, the
routing matrix $P := [p_{i,j}]$ is sub-stochastic. Note that, we also
allow feedback of jobs to the same node; i.e., $p_{i,i} \geq 0$. Let
$\phi_l^i : \Omega \rightarrow \{1,\ldots,K\}, \,\, \forall i \in\sK \text{ and } \forall l \in \mathbb{N}$, be a measurable
function such that $\phi_l^i = j$ implies that the $l$th job at node
$i$ will be routed to node $j$ and $\mathbb E[\mathbf 1_{\{\phi_l^i = j\}}] = p_{i,j}$.% For brevity, let $\Gamma_i^{k} = e_{\phi_{i}^k}$, where $e_i$ is the $i$th unit vector.
~Define the random vector $R_l(m) := \sum_{i=1}^m e_{\phi_i^l}$, where
$e_i$ is the $i$th $K$-dimensional unit vector and the $k$th component
of $R_l(m)$, denoted $R_l^k(m)$, represents the number of departures from node $l$ to node $k$ out of $m$ departures from that node. Then, $\mathbf R(m) := (R_1(m), \ldots, R_K(m))$ is a $K \times K$ matrix whose columns are the routing vectors from the nodes in the network.
\begin{assumption} \label{thm:routing-limits}
The stationary routing process $\{\mathbf R(m),~m \geq 1\}$ satisfies the
following functional limits:

\noindent (i) 
\(
n^{-1}\mathbf R(n e) \to  \mathbf P~e 
\)
in $(\sC^{K \times K}, U)$ a.s. as $n \to \infty$, where $e: [0,\infty) \to
[0,\infty)$ is the identity function, and

\noindent (ii) 
\(
\hat{\mathbf R}_n := n^{-1/2}\left( \mathbf R(n e) -
  n \mathbf P e\right) \Rightarrow \hat {\mathbf R},
\)
in $(\sC^{K \times K},U)$ as $n \to \infty$, where $\hat{\mathbf R} = [W_{i,j}]$ and
$W_{i,j}$ are independent Brownian motion processes with mean zero and
diffusion coefficient $p_{i,j}(1-p_{i,j})$.
\end{assumption}

As a direct consequence of Assumption~\ref{thm:routing-limits} we have
the following corollary, which will prove useful in our analysis of
the network state process in the next section.

\begin{corollary} \label{cor:routing-limits}
The routing process $\mathbf R$ also satisfies the following fCLT:
\begin{equation}
  \label{eq:2}
  \hat{\mathbf R}_n^{T} \mathbf 1 \Rightarrow \hat{\mathbf R}^{T}\mathbf 1 = \tilde{\mathbf{W}} ~\text{in}~(C^K,U)~\text{as}~n \to 
\infty,
\end{equation}
where $\mathbf 1 = (1,\ldots, 1)$ is a $K$-dimensional vector of one's
and 
\[
\tilde{\mathbf W} = \left(\sum_{k=1}^K W_{1,k}, \ldots, \sum_{k=1}^K
W_{K,k} \right)
\]
 is a $K$-dimensional Brownian motion with mean zero and
diagonal covariance matrix with entries
\[
\left(\sum_{k=1}^K p_{1,k}(1 - p_{1,k}), \ldots, \sum_{k=1}^K p_{K,k}(1-p_{K,k}) \right).
\]
\end{corollary}

Finally, we claim the following joint convergence result that
summarizes and generalizes the convergence results in the
afore-mentioned theorems.
\begin{proposition} \label{prop:primitive-limits}
Assume that for each $n \geq 1$, $\mathbf A_n$, $\mathbf S_n$ and $\mathbf R(n)$ are mutually independent. Then,\\
\noindent (i) 
\(
n^{-1}\left( \mathbf A_n, \mathbf S_n, \mathbf R(n e)\right) \to \left( \mathbf F, \mathbf M, \mathbf P{'} e \right) 
\)
in $(\sC^J\times\sC^K\times \sC^{K\times K}, U)$ a.s. as $n \to \infty$, and\\
\noindent (ii) 
\(
\left( \hat{\mathbf A}_n, \hat{\mathbf S}_n, \hat{\mathbf R}_n\right) \Rightarrow \left( \mathbf W^0 \circ \mathbf F, \mathbf W \circ \mathbf M, \hat{\mathbf R}\right)
\)
in $(\sC^J\times\sC^K\times \sC^{K\times K}, U)$ as $n \to \infty$.
\end{proposition}
The joint convergence follows from the assumed independence of the
pre-limit random variables, and is straightforward to establish under
the uniform convergence criterion.

\section{Functional Strong Law of Large Numbers} \label{sec:fluid-networks}
Let $Q_{n,k}(t) = E_{n,k}(t) - D_{n,k}(t)$ be the queue length sample path at node
$k$, where $E_{n,k}(t)$ is the total number of jobs arriving at node $k$ in the interval $[0,t]$
and $D_{n,k}$ is the cumulative departure process. We assume that the
server is non-idling implying that $D_{n,k}(t)= S_{n,k}(B_{n,k}(t))$, where $B_{n,k}(t) := \int_0^t \mathbf 1_{\{Q_{n,k}(s) > 0 \}} ds$ is the total busy time of the server. Therefore, the queue length process is
\begin{eqnarray} \label{queue-length-path}
Q_{n,k}(t) := A_{n,k}(t) + \sum_{l=1}^K R_{l}^k(S_{n,l}(B_{n,l}(t))) - S_{n,k}(B_{n,k}(t)).
\end{eqnarray}

% Recall the queue length process sample path (of node $k \in \sK$)
% defined in \eqref{queue-length-path}. 
The $K$-dimensional multivariate
stochastic process \(\mathbf Q_n := (Q_{n,1}, \ldots, Q_{n,K}) \in \sD^K\) represents the
network state. Our first result establishes a fluid limit
approximation to a rescaled version of $\mathbf Q_n$ by establishing a functional strong law
of large number result as the exogenous arrival population size $n$
scales to infinity. Consider the queue length in the $k$th node,
$Q_k$. Centering each term on the right hand side by the corresponding
fluid limits (and subtracting those terms), and introducing the term
$\int_0^t \mu_{n,k}(s) ds$, we obtain $n^{-1} Q_{n,k}(t)$
\begin{eqnarray}
\nonumber
& &\begin{split}
=  & \left( \frac{1}{n}A_{n,k}(t) - F_k(t) \right) + \left(\frac{1}{n}
\sum_{l=1}^K \bigg [ R_{l}^k (S_{n,l}(B_{n,l}(t))) -  p_{l,k}
S_{n,l}(B_{n,l}(t))\bigg ] \right) \\ &- \left(\frac{S_{n,k}(B_{n,k}(t))}{n} - \int_0^{B_{n,k}(t)} \mu_{n,k}(s)
  ds \right)\\& + \left( F_k(t) -  \int_0^{B_{n,k}(t)} \mu_{n,k}(s)
  ds + \frac{1}{n}
\sum_{l=1}^K  p_{l,k} S_{n,l}(B_{n,l}(t)) \right) 
\end{split}\\
\label{queue-length-fluid scale}
& &\begin{split}  = & \left(\frac{1}{n}A_{n,k}(t) - F_k(t) \right)+
  \left(\frac{1}{n} \sum_{l=1}^K \bigg [ R_{l}^k (S_{n,l}(B_{n,l}(t)))
    -  p_{l,k} S_{n,l}(B_{n,l}(t))\bigg ] \right) \\
&- \left(\frac{S_{n,k}(B_{n,k}(t))}{n} - \int_0^{B_{n,k}(t)} \mu_{n,k}(s)
  ds \right)\\
&+ \left(F_k(t) - \int_0^t \mu_{n,k}(s) ds \right) + (1-p_{k,k}) \int_{B_{n,k}(t)}^{t}\mu_{n,k}(s) ds \\
& + \left(\frac{1}{n} \sum_{l=1}^K p_{l,k} \bigg [ S_{n,l}(B_{n,l}(t)) - n
  \int_0^{B_{n,l}(t)} \mu_{n,l}(s) ds \bigg ] \right) \\
&+ \sum_{l=1}^K p_{l,k} \left( \int_0^t \mu_{n,l}(s) ds \right)
- \sum_{l \ne k} p_{l,k}
\int_{B_{n,l}(t)}^t \mu_{n,l}(s) ds. \\
\end{split}
\end{eqnarray}
Note that we used the fact that $B_{n,k}(t) \leq t$ so that
\[
\int_0^t \mu_{n,k}(s) ds = \int_0^{B_{n,k}(t)} \mu_{n,k}(s) ds +
\int_{B_{n,k}(t)}^t \mu_{n,k}(s) ds.
\]
Recall too that $I_{n,k}(t) := t
- B_{n,k}(t) = \int_{0}^t
\mathbf 1_{\{Q_{n,k}(s) = 0\}} ds$ is
the idle time process, which measures the amount of time in
$[0,t]$ that the node is not serving jobs (i.e., the queue is
empty). Now, $n^{-1}  Q_{n,k}$ can be decomposed as the sum of two
processes, $\bar X_{n,k}$ and $\bar Y_{n,k}$, where
\begin{eqnarray}
\label{X-fluid-scale}
&& \begin{split} \bar{X}_{n,k}(t)= & \left(\frac{1}{n}A_{n,k}(t) -
    F_k(t) \right)\\ &+ \left(\frac{1}{n} \sum_{l=1}^K \bigg [ R_{l}^k (
    S_{n,l}(B_{n,l}(t))) - p_{l,k} S_{n,l}(B_{n,l}(t))\bigg ]
  \right)\\
&- \left(\frac{S_{n,k}(B_{n,k}(t))}{n} - \int_0^{B_{n,k}(t)}
  \mu_{n,k}(s) ds \right)\\
&+ \left(F_k(t) -
  \left( \int_0^t \mu_{n,k}(s) ds \right) \right) + \sum_{l=1}^K p_{l,k} \left( \int_0^t \mu_{n,l}(s) ds \right)\\
&+ \left( \frac{1}{n} \sum_{l=1}^K p_{l,k} \bigg [S_{n,l}(B_{n,l}(t))-
  n\int_0^{B_{n,l}(t)} \mu_{n,l}(s) ds \bigg ] \right), \text{ and }
\end{split}\\
\label{Y-fluid-scale}
&&\bar{Y}_{n,k}(t)= (1-p_{k,k}) \int_{B_{n,k}(t)}^t \mu_{n,k}(s) ds - \sum_{l \ne k} p_{l,k} \int_{B_{n,l}(t)}^t \mu_{n,l}(s) ds.
\end{eqnarray}

While this expression appears formidable, the analysis is
simplified significantly by the fact that $\bar{\mathbf {Q}}_n :=
n^{-1}(Q_{n,1}, \ldots, Q_{n,k})$ and $\bar{\mathbf Y}_n := (\bar Y_{n,1}, \ldots, \bar
Y_{n,k})$ are solutions to the $K$-dimensional Skorokhod/oblique
reflection problem. First, we recall the definition of the oblique reflection problem.

\begin{theorem} \label{thm:ORT}[Oblique Reflection Problem] 
Let $\mathbf V$ be a $K \times K$ $M$-matrix\footnote{An $M$-matrix is
  a square matrix with spectral radius less than one.}, also known as
the reflection matrix. Then, for
every $x \in \mathcal{D}_0^K := \{x \in \mathcal{D}^K : x(0) \geq
0\}$, there exists a unique tuple of functions $(y,z)$ in
$\mathcal{D}^{K} \times \mathcal D^K$ satisfying
\begin{eqnarray}
\nonumber
z &=& x + \mathbf V y \geq 0,\\
\label{eq:orm}
dy &\geq& 0 \text{ and } y(0) = 0,\\
\nonumber
z_j dy_j &=& 0, \quad j = 1, \ldots, K.
\end{eqnarray}
The process $(z,y) := (\Phi(x),\Psi(x))$ is the so-called oblique reflection map, where $\Phi(x) = x + \mathbf V \Psi(x)$.
\end{theorem}
Note that, in general, if $\mathbf G$ is a nonnegative $\mathbf
M$-matrix then so is $\mathbf V = \mathbf I- \mathbf G^T$ (Lemma
7.1 of \cite{ChYa01}). The following lemma shows that the queue length
satisfies the Oblique Reflection Mapping.

\begin{lemma}
\label{lem:queue-ORT}
Consider $\mathbf {\bar{X}}_n(t) = (\bar{X}_{n,1}(t), \ldots,
\bar{X}_{n,k}(t)) \in \mathcal{D}_0^K$, where $\bar{X}_{n,k}(t) \,\, k
\in \{1,\ldots,K\}$ is defined in \eqref{X-fluid-scale},
$\mathbf{\bar{Q}}_n \in \mathcal{D}^K$ and $\mathbf{\bar{Y}}_n \in \mathcal{D}_0^K$. Then,
\(
(\mathbf{\bar{Q}}_n, \mathbf{\bar{Y}}_n) = (\Phi(\mathbf{\bar{X}}_n), \Psi(\mathbf{\bar{X}}_n)),
\)
with reflection matrix $\mathbf V = \mathbf I - \mathbf P^T$.
\end{lemma}

Next, we establish a functional strong law of
large numbers result for \eqref{X-fluid-scale}, which will subsequently be
used in Theorem~\ref{thm:queue-fluid} for the queue length approximation.
\begin{lemma} \label{lem:X-fluid}
  The fluid-scaled netput process $\bar{\mathbf X}_n$ converges
  to a deterministic limit as $n \to \infty$:
\begin{eqnarray*}
\mathbf{\bar{X}}_n(t) \rightarrow \bar{\mathbf X}(t) := (\bar{X}_1(t),
  \ldots, \bar{X}_K(t)) ~\text{in}~(\sC^K,U)~\text{a.s.},
\end{eqnarray*}
where, 
\begin{eqnarray}
\label{X-fluid}
\bar{X}_k(t) = F_k(t) - \int_0^t \mu_k(s)ds + \sum_{l=1}^K p_{l,k} \int_0^t \mu_l(s) ds.
\end{eqnarray}
\end{lemma}

We can now establish the functional strong law of large numbers limit
for the queue length process. The proof essentially follows from the
continuity of the oblique reflection map $(\Phi(\cdot), \Psi(\cdot))$.
\begin{theorem} \label{thm:queue-fluid}
Let $\mathbf{\bar{X}}_n(t)$ and $\mathbf{\bar{X}}(t)$ be as defined in \eqref{X-fluid-scale} and \eqref{X-fluid} respectively. Then, $(\mathbf{\bar{Q}}_n(t), \bar{\mathbf Y}_n(t))$ satisfy Theorem \ref{thm:ORT} and, as $n \rightarrow \infty$,
\[
(\bar{\mathbf Q}_n(t), \bar{\mathbf Y}_n(t)) \rightarrow
(\Phi(\bar{\mathbf X}(t)), \Psi(\bar{\mathbf
  X}(t)))~\text{in}~(\sC^K\times\sC^K)~\text{a.s.}~\forall t \in [0,\infty).
\]
\end{theorem}
\Proof
It follows by Lemma \ref{lem:queue-ORT} that $(\bar{\mathbf Q}_n(t), \bar{\mathbf Y}_n(t))$ satisfy the oblique reflection mapping theorem. Therefore,
\(
(\bar{\mathbf Q}_n(t), \bar{\mathbf Y}_n(t)) \equiv (\Phi(\bar{\mathbf X}_n(t)), \Psi(\bar{\mathbf X}_n(t))).
\)
Now, the reflection regulator map, $\Psi(\cdot)$, is Lipschitz continuous under the uniform metric (Theorem 7.2, \cite{ChYa01}). By the Continuous Mapping Theorem and Lemma \ref{lem:X-fluid} it follows that,
\[
(\Phi(\bar{\mathbf X}_n(t)), \Psi(\bar{\mathbf X}_n(t))) \rightarrow (\Phi(\bar{\mathbf X}(t)), \Psi(\bar{\mathbf X}(t))) \,\, u.o.c. \,\, a.s. \text{ as } n \rightarrow \infty, \,\, \forall t \in [0,\infty).
\]
\EndProof

\Proof [Proof of Lemma~\ref{lem:queue-ORT}]
First, by definition we have $\mathbf{\bar{Q}}_n = \mathbf{\bar{X}}_n
+ (I - \mathbf P^{T}) \mathbf{\bar{Y}}_n$. Note that $\mathbf P$ is a
non-negative (sub-stochastic) matrix with spectral radius less than
unity and, therefore, an $M$-matrix, implying that $\mathbf I- \mathbf
P^{T}$ is also an $M$-matrix. Once again by definition $Q_{n,k}$ and
$Y_{n,k}$ satisfy the conditions in \eqref{eq:orm} for all $k  \in
\sK$. Thus, Theorem \ref{thm:ORT} is satisfied and the lemma is proved.
\EndProof

\Proof [Proof of Lemma~\ref{lem:X-fluid}]
The result follows by an application of part (i) of Proposition
\ref{prop:primitive-limits} to \eqref{X-fluid-scale}. Noting that
$B_{n,k}(t) \leq t$, the random time change theorem (Theorem 5.5,
\cite{ChYa01}) and Assumption~\ref{thm:service-limits}(i) together imply that,
\[
\frac{1}{n}S_{n,k}(B_{n,k}(t)) - \int_0^{B_{n,k}(t)}
  \mu_{n,k}(s) ds \rightarrow 0~\text{in}~(\sC,U)~\text{a.s. as}~ n \to \infty \,\, \forall t \in [0,\infty).
\]
Similarly, applying the random time change theorem along with
Assumption~\ref{thm:routing-limits} (i) and Assumption~\ref{thm:service-limits}(i)
we obtain
\[
\frac{1}{n} \left(R_{l}^k( S_{n,k}(B_{n,k}(t))) - p_{l,k}
  S_{n,k}(B_{n,k}(t)) \right) \rightarrow
0~\text{in}~(\sC,U)
\]
$~\text{a.s. as}~ n \rightarrow \infty \,\, \forall t \in [0,\infty).$
Applying these results to \eqref{X-fluid-scale} it follows that
$\bar{X}_{n,k}(t) \rightarrow \bar{X}_k(t)~\text{in}~(\sC,U)~\text{a.s. as}n
\rightarrow \infty$. The joint convergence follows automatically from
these results and Proposition~\ref{prop:primitive-limits}.
\EndProof

Note that neither Theorem \ref{thm:ORT} nor Theorem
\ref{thm:queue-fluid} provide an explicit functional form for the
reflection regulator $\Psi(\cdot)$. It can be shown (see \cite[Chapter 7]{ChYa01}) that the regulator
map is the unique fixed point, $y^* \in \mathcal{D}^K$, of the map
$\pi(x,y)(t) := \sup_{0 \leq s \leq t} [-x(s) + \mathbf G y(s)]^+ \,\, \forall
t \in [0,\infty)$, where $\mathbf G$ is an $M$-matrix. Extracting a closed form expression
for $y^*$ is not straightforward, barring a few special cases. The following corollary shows that the reflection map and fluid limit of the queue length process for a parallel node queueing network is particularly simple and an obvious generalization of that of a single queue.

\begin{corollary} \label{cor:queue-fluid}
Consider a $K$-node parallel queueing network. The fluid limit to the
queue length and cumulative idleness processes are
\[
(\bar{\mathbf Q}, \bar{\mathbf Y}) = (\Phi(\bar{\mathbf X},
\Psi(\bar{\mathbf X}))) \in \sC^K \times \sC^K,
\]
where $\bar X = (X_1, \ldots, X_K)$, $\Psi(\bar{\mathbf X}(t)) =
\sup_{0 \leq s \leq t} [-\bar{\mathbf X}(s)]^+$ and $\Phi(\bar{\mathbf
X}) = \bar{\mathbf X} + \Psi(\bar{\mathbf X})$.
\end{corollary}
\Proof
  Note that for a parallel queueing network $\mathbf P =
  0$. Therefore, the fixed point of the map $\pi(\cdot,\cdot)$ is
  simply $\sup_{0\leq s \leq t} [-x(s)]^+$. It follows that the regulator map of the fluid scaled queue length process is
\(
\Psi(\bar{\mathbf X}_n(t)) = \sup_{0 \leq s \leq t} [-\bar{\mathbf X}_n(s)]^+.
\)
It follows by Theorem \ref{thm:queue-fluid} that \[\Psi(\bar{\mathbf X}_n(t)) \rightarrow \sup_{0 \leq s \leq t} [-\bar{\mathbf X}(s)]^+\] and \[\Phi(\bar{\mathbf X}_n(t)) \rightarrow \bar{\mathbf X}(t) + \Psi(\bar{\mathbf X}(t))~\text{in}~(\sC^K,U)~\text{a.s.}\] as $n \rightarrow \infty$. 
\EndProof

A slightly more complicated example would be a series queueing
network.  Corollary~\ref{cor:tandem-fluid} establishes the fluid limit
to the network state of a two queue tandem network, when a large, but finite, number $n$ of users
arrive at queue 1 over a finite time horizon $[-T_0,T]$. This result
can be rather straightforwardly extended to a network of more than two
queues. Let  \(
\mathbf P = \begin{pmatrix} 
0 & 0 \\
1 & 0
\end{pmatrix},
\) be the matrix of Markov routing probabilities.

\begin{corollary} \label{cor:tandem-fluid}
Consider a tandem queueing network 
and recall that $\mathbf V = \mathbf I - \mathbf P^T$. Let $\mathbf F = F_1$ be the arrival epoch distribution with support $[-T_0,T]$ where
$T_0,T > 0$, and assume that $\m_1$ and $\m_2$ are the fixed service
rates. Then, the (joint) fluid limit to the queue length and
cumulative idleness processes is
\(
(\bar{\mathbf Q}, \bar{\mathbf Y}) = (\Phi(\mathbf{\bar X}), \Psi(\mathbf{\bar X})) \in \sC^K\times\sC^K,
\)
where
\[
\mathbf {\bar X} := (X_1, X_2) = ((F_1 - \m_1 e), (\m_1 - \m_2) e), 
\]
and 
~\[
\Psi(\bar{\mathbf X}) = (Y_1, Y_2)
\]
with $Y_1(t) = \sup_{0 \leq s \leq t}(-X_1(s))_+$ and $Y_2(t) = \sup_{0\leq s
  \leq t} (-X_2(s) + Y_1(s))_+ = \sup_{0 \leq s \leq t} [-X_2(s) + \sup_{0\leq r \leq s}
(-X_1(r))_+]_+$, and $\Phi(\mathbf{\bar X})(t) =\bar{\mathbf X} +
\mathbf V \Psi(\bar{\mathbf X}) = (X_1 + Y_1, ~X_2 + Y_2 - Y_1)$.
\end{corollary}
The proof is straightforward by substitution and we omit it. Note that
the queue length fluid limit to the downstream queue appears quite
complicated: $\bar Q_2 = X_2 + Y_2 - Y_1$ where $Y_2(t) = \sup_{0 \leq
s \leq t} (- X_2(s) + Y_1(s))_+$. By substituting in the expression
for $X_2$ we have 
\begin{eqnarray*}
\bar Q_2 &=& (\m_1 - \m_2)e + F_1 - F_1 - Y_1 + Y_2\\
&=& (F_1 - \bar Q_1 -\m_2 e) + Y_2.
\end{eqnarray*}
Note that $F_1 - \bar Q_1$ is just the cumulative fluid departure
function from the upstream queue, which is precisely the input to the
downstream queue. 

Next, suppose the service process is stationary such that $\mu_k(t) =
\mu_k$ for all $t \geq 0$ and $k \in \sK$. Then, the busy time process
satisfies the following theorem.
\begin{theorem}
\label{thm:busy-time-fluid}
Let $\bar{\mathbf B}_n(t) = (B_{n,1}(t), \ldots, B_{n,k}(t))$. Then, as $n \rightarrow \infty$,
\begin{eqnarray} \label{busy-time-fluid}
\bar{\mathbf B}_n \rightarrow e \mathbf{1} - \mathbf M \Psi(\bar{\mathbf X})~\text{in}~(\sC^K,U) ~\text{a.s.},
\end{eqnarray}
where, $\mathbf M = diag(1/\mu_1, \ldots, 1/\mu_K)$.
\end{theorem}
\Proof
By definition 
\(
\bar{\mathbf B}_n(t) = t \mathbf{1} - \bar{\mathbf I}_n(t),
\)
where $\bar{\mathbf I}_n(t) = (I_{n,1}(t), \ldots, I_{n,k}(t))^T$. Recalling the definition of the process $\mathbf {\bar Y}_n(t)$ it is straightforward to see that $\mathbf {\bar I}_n(t) = (\mathbf I - \mathbf P^T)^{-1} \mathbf{\bar Y}_n(t)$ for all $t \geq 0$. Therefore,
\(
\bar{B}_n(t) = \underbar{t} - (\mathbf I - \mathbf P^T)^{-1} \bar{\mathbf Y}_n(t).
\)
~Theorem \ref{thm:queue-fluid} implies that, as $n \rightarrow \infty$,
\(
\bar{\mathbf B}_n(t) \rightarrow \underbar{t} - \Psi(\bar{\mathbf X}(t)) ~\text{in}~(\sC^K,U)~\text{a.s.}~\forall t \in [0,\infty).
\)
\EndProof

The following corollary establishes the fluid busy time process for the parallel queue case. The proof follows that of Corollary \ref{cor:queue-fluid} and we omit it.
\begin{corollary}
\label{cor:busy-time-fluid}
Consider a $K$-node parallel queueing network. Then,
\[
\bar{\mathbf B}_n \rightarrow e \mathbf{1} - (\mathbf I - \mathbf
P^T)^{-1} \sup_{0 \leq s \leq \cdot} [-\bar{\mathbf X}(s)]^+
 ~\text{in}~(\sC^K,U)~\text{a.s.}~n \to \infty.
\]
\end{corollary}
In the stationary
case we considered here, the busyness time-scale is effectively fixed
by the service rate through the matrix $\mathbf M$. On the other hand, if the
service processes are non-stationary this time-scale \textit{itself} is
time-varying. Thus, computing the
busy time (or equivalently the idle time) process when the service
process is non-stationary is complicated. Note that the function $\bar{
\mathbf{Y}}$ represents the number of ``blanks'' or the amount of
unused capacity in the network at each point in time, providing an
indication of whether a particular queue in the network is busy or
not.

Note that the population acceleration scale we use in the current
analysis ensures that (in the limit) the amount of time each user
spends in service is infinitesimally small. The `behavior'
of the queue state under population acceleration scaling is akin
to the conventional heavy-traffic scaling introduced in \cite{Re1984} for
stationary single class queueing networks. The corresponding diffusion heavy-traffic scaling
identifies the critical time-scale of the
stationary queueing network. The population acceleration scaling
differs from the conventional heavy-traffic scaling by the fact that
the fluid limit process (in general) is non-linear in nature. This
implies that queues in the network can enter idle and busy periods, and arriving jobs will only
face delays in the latter time intervals. We should anticipate that the
critical time-scale of the
queue state in the diffusion scale should itself change
depending on whether the queue is busy or idle, leading to a
non-stationary diffusion approximation. Indeed, this is
precisely what is implied by the results in the next section.

\section{Functional Central Limit Theorems} \label{sec:diffusion-networks}
We now consider the second order refinement to the fluid
limit by establishing a functional central limit theorem (FCLT) satisfied by
the queue length state process. We show, in particular, that the 
FCLT is a reflected diffusion, where the diffusion process $\hat{\mathbf{X}}$ is a
function of the multi-dimensional Brownian bridge process as in Definition~\ref{def:bridge}. Unlike the heavy traffic limits for
generalized Jackson networks (see \cite[Chapter 7]{ChYa01}
\cite{Re1984}), the diffusion is \textit{not} reflected through
the oblique reflection map (see \cite[Definition 7.1]{ChYa01}). The non-homogeneous
traffic and non-stationary service processes induce a time-varying critical
time-scale under the population acceleration scaling.
we show that this time-varying critical time-scale manifests as a
time-varying reflection boundary in transitory queueing networks. To be precise, the reflection
regulator for the queue length diffusion  is the directional derivative of the oblique reflection of $\bar{\mathbf{X}}$ (from
Lemma \ref{lem:X-fluid}) in the direction of the diffusion limit $\hat{
\mathbf{X}}$ to the netput process. A similar result was observed in the case of
a single $\D_{(i)}/GI/1$ transitory queue in \cite{HoJaWa2012}.%  In that case,
% the directional derivative reflection map was explicitly
% characterized by appealing to the results in \cite[Chapter
% 9]{Wh01b}. On the other hand, the results in \cite{MaRa10}
% characterize the directional derivative of the multidimensional
% oblique reflection map.

Recall that $\mathbf V$ is a $K \times K$ $M$-matrix and $\mathbf P^T =
\mathbf I- \mathbf V$.
Let $x \in \sC_0$ then, under the hypothesis of Theorem
\ref{thm:ORT}, there exists a unique \textit{oblique reflection map} $(z,y) := (\Phi(x),\Psi(x)) \in \sC \times \sC$ such that $z = x
+ \mathbf V y$, $y_j$ is non-decreasing and $y_j$
grows only when $z_j$ is zero (for all $j = 1,\ldots, K$). The
directional derivative of the oblique reflection of $x$ in the
direction of the process $\chi \in \sC$ is defined as follows (see
\cite{MaRa10} as well):

\begin{definition} [Directional Derivative Reflection Map] \label{def:dir-der}
Given $(x,\chi) \in \sC^K_0 \times \sC^K$ and $M$-matrix $\mathbf V$, the directional derivative of the oblique reflection map $\Phi(x) = x
+ \mathbf V \Psi(x)$ in the direction of $\chi$ is the pointwise limit of
\[
\mathbf \Delta_\chi^n(x) := \sqrt{n} \left( \Phi\left( \frac{\chi}{\sqrt n} + x \right) - \Phi(x)
\right) \in \sC ~~ n \geq 1
\]
as $n \to \infty$.
\end{definition} 

Theorem 1.1 (ii) of \cite{MaRa10} identifies the limit process, which we
state as a lemma for completeness. Here,

\begin{lemma} \label{lem:dir-der-ORT}
If $(x,\chi) \in \sC^K_0 \times \sC^K$ then the directional derivative
limit $\mathbf \D_{\chi}(x)$ exists and convergence in Definition
\ref{def:dir-der} is uniformly on compact subsets of continuity points
of the limit $\mathbf \D_\chi(x)$. Further, if $(z,y)$ solve the
oblique reflection problem for $x$ then
\[
\mathbf \D_\chi(x) = \chi + \mathbf V \g (x,\chi),
\]
where $\mathbf V = \mathbf I - \mathbf P^T$ $\g := \g(x,\chi)$ lies in $\sD^K_{usc}$ and is the unique
solution to the system of equations
\[
\g^i(t) = \begin{cases}
\sup_{s \in \nabla_t^i} [-\chi^i(s) + [\mathbf P \g]^i(s)]_+ &  t \in [0, t_u^i], \\
\sup_{s \in \nabla_t^i} [-\chi^i(s) + [\mathbf P \g]^i(s)] & t > t_u^i,
\end{cases}
\]
for $i = 1, \ldots, K$, where 
\(
\nabla_t^i := \{s \in [0,t] | z^i(s) = 0 \text{ and } y^i(s) = y^i(t)\},
\)
and 
\(
t_u^i := \inf \{t \geq 0: y^i(t) > 0\}.
\)
\end{lemma}

Consider the second order refinement to the netput process,
\[
\hat{\mathbf X}_n := n^{-1/2} \left(  {\mathbf X}_n- n \bar{\mathbf X} \right) \in \sD^K.
\]
Using Proposition \ref{prop:primitive-limits}, and the fact that the
limit processes have sample paths in $\sC^K$, the following Lemma is
straightforward to establish.%  We abuse notation slightly and denote
% composition of two vector-valued functions as $x \circ y = (x_1 \circ
% y_1, \ldots, x_K \circ y_K)$.

\begin{lemma} \label{lem:X-hat}
The diffusion-scaled netput process satisfies,
\[
\hat{\mathbf X}_n \Rightarrow \hat{\mathbf X} ~\text{in}~ (\sC^K,U) \text{ as } n \to \infty,
\]
where $\hat X_k := W^0_k \circ F_k - W_k \circ \int_0^t \mu_k(s) ds +
\mathbf 1^T (\hat{\mathbf R}_k \circ\mathbf M)$,
$\hat{\mathbf R}_k$ is the $k$th row of the matrix valued process
$\hat{\mathbf R}$ defined in part (ii) of
Assumption~\ref{thm:routing-limits}, and $\mathbf M$ is defined in \eqref{eq:1}.
\end{lemma}
The proof of the lemma follows from an application of part (ii)
of Proposition~\ref{prop:primitive-limits} and using the fact that the
addition operator is a continuous map under the uniform metric. We
omit it for brevity.

Now, The diffusion scale queue length process is 
\[
\hat{\mathbf Q}_n := n^{-1/2} \left( {\mathbf Q}_n - n \mathbf {\bar Q} \right) \in \sD^K.
\]
Recall, from Lemma \ref{lem:queue-ORT}, that $\mathbf {\bar{Q}}_n =
\mathbf{\bar X}_n + \mathbf V \Psi(\mathbf{\bar X}_n)$ and, from
Theorem \ref{thm:queue-fluid}, that $\mathbf{\bar Q} = \mathbf{\bar X}
+ \mathbf V \Psi(\mathbf{\bar X})$. It follows that 
\begin{eqnarray*}
\mathbf{\hat Q}_n &=& \sqrt{n} \left( \mathbf{\bar X}_n + \mathbf V
  \Psi(\mathbf{\bar X}_n) -  \mathbf{\bar X}
- \mathbf V \Psi(\mathbf{\bar X})\right)\\
&=&
    \begin{split}
      \hat{\mathbf X}_n &+ \mathbf V \sqrt n \left( \Psi \left(
          \frac{\hat{\mathbf X}_n}{\sqrt n} + \bar{\mathbf X} \right)
        - \Psi \left(\bar{\mathbf X} \right) \right)\\ &+ \mathbf V \sqrt
      n \left( \Psi(\bar{\mathbf X}_n) - \Psi
        \left(\frac{\hat{\mathbf X}_n}{\sqrt n} + \bar{\mathbf X}
        \right) \right)
    \end{split}
\\
&=& \D^n_{\hat{\mathbf X}_n}\left( \bar{\mathbf X} \right) +  \mathbf V \sqrt n \left(
    \Psi(\bar{\mathbf X}_n) - \Psi \left(\frac{\hat{\mathbf
    X}_n}{\sqrt n} + \bar{\mathbf X} \right) \right).
\end{eqnarray*}

% Now, suppose that
% $\mathbf {\hat X} \in \sC^K$. 
Our
next result shows that $\mathbf \D^n_{\mathbf{\hat X}_n}(\mathbf{\bar
  X})$ is asymptotically equal to $\mathbf \D^n_{\mathbf{\hat X}}(\mathbf{\bar X})$.

\begin{lemma} \label{lem:uni-dir-der}
Let $\mathbf \D^n_{\mathbf{\hat X}}(\mathbf{\bar X})$ and $\mathbf
\D^n_{\mathbf{\hat X}_n}(\mathbf{\bar X})$ be defined as in Definition
\ref{def:dir-der}. Then,
\[
\mathbf \D^n_{\mathbf{\hat X}_n}(\mathbf{\bar X}) - \mathbf
\D^n_{\mathbf{\hat X}}(\mathbf{\bar X}) \to 0 ~\text{in}~(\sC^K,U)~\text{a.s. as}~ n \to \infty.
\]
\end{lemma}

% In this case Lemma \ref{lem:dir-der-ORT}
% implies that  $\mathbf \D^n_{\mathbf{\hat X}}(\mathbf{\bar X})$
% converges  $\mathbf \D_{\mathbf{\hat X}}(\mathbf{\bar X})$, the
% directional derivative limit.
Lemma~\ref{lem:uni-dir-der} implies it suffices to
consider the process
\begin{equation} \label{eq:new-diff-queue}
\hat{\mathbf Q}_n \equiv \D^n_{\hat{\mathbf X}}\left( \bar{\mathbf X} \right) +  \mathbf V \sqrt n \left(
    \Psi(\bar{\mathbf X}_n) - \Psi \left(\frac{\hat{\mathbf
    X}}{\sqrt n} + \bar{\mathbf X} \right) \right)
\end{equation}
(where, by an abuse of notation, we call this process $\hat{\mathbf
  Q}_n$ as well). Now, if we
show that 
\[
\sqrt n \left(
    \Psi(\bar{\mathbf X}_n) - \Psi \left(\frac{\hat{\mathbf
    X}}{\sqrt n} + \bar{\mathbf X} \right) \right) \to 0~\text{in}~(\sC^K,U)
\]
a.s. as $n \to \infty$, then Lemma~\ref{lem:dir-der-ORT} implies that $\hat{\mathbf Q}_n$
converges to the process
$\mathbf \D_{\mathbf{\hat X}}(\mathbf{\bar X})$ pointwise in the large
population limit. The following lemma establishes the required result
under general conditions.

\begin{lemma} \label{lem:extra-terms}
  Let $x_n, x \in \sD^K$ be stochastic processes that satisfy $\| \sqrt n (x_n -x)\| \to
  \chi$ a.s. as $n \to \infty$. Then, 
  \begin{equation}
    \label{eq:extra-terms}
    \left \| \sqrt n \left( \Psi(x_n) - \Psi \left( \frac{\chi}{\sqrt
            n} + x\right) \right) \right\| \to 0 ~\text{a.s. as}~n \to \infty,
  \end{equation}
  where $\chi \in \sC^K$.
\end{lemma}

The main result of this section follows as a consequence of these lemma's.
\begin{theorem} \label{thm:pt-wise}
 For any fixed $t \in [0,\infty)$, as $n \to \infty$
  \begin{equation}
    \label{eq:pt-wise}
    \hat{\mathbf Q}_n(t) \Rightarrow \hat{\mathbf Q}(t) =
    \D_{\hat{\mathbf X}}(\bar{\mathbf X})(t),
  \end{equation}
  where $\D_{\hat{\mathbf X}}(\bar{\mathbf X})(t) = \hat{\mathbf
    X}(t) + \mathbf V \gamma(\bar{\mathbf X}, \hat{\mathbf X})(t)$.
\end{theorem}
\Proof
  First, using the Skorokhod representation theorem \cite{Bi68}, it
  follows from Lemma~\ref{lem:X-hat} that there exist versions
  of the stochastic processes $\left \{ \hat{\mathbf X}_n \right\}$
  and $\hat{\mathbf X}$, referred to using the same notation, such
  that
  \[
  \hat{\mathbf X}_n\to \hat{\mathbf X} ~\text{in}~ (\sC^K,U) ~\text{a.s. as}~ n \to \infty.
  \]
It follows that $\bar{\mathbf X}_n = \bar{\mathbf X} + (\sqrt n)^{-1}
\hat{\mathbf X} + o(1)$ a.s. Lemma~\ref{lem:extra-terms} implies that 
\[
\left \| \sqrt n \left(
    \Psi(\bar{\mathbf X}_n) - \Psi \left(\frac{\hat{\mathbf
    X}}{\sqrt n} + \bar{\mathbf X} \right) \right) \right \| \to 0
\]
a.s. as $n \to \infty$. Next, using Lemma~\ref{lem:uni-dir-der} and
Lemma~\ref{lem:dir-der-ORT}, it follows that $\hat{\mathbf Q}_n(t) \to
\D_{\hat{\mathbf X}}(\bar{\mathbf X})(t) $ a.s. as $n \to
\infty$ for any fixed $t \in [0,\infty)$, which in turn implies weak convergence of the stochastic
processes thus proving the desired result.
\EndProof

\Proof[Proof of Lemma~\ref{lem:uni-dir-der}]
First, recall that
\[
\mathbf \D^n_{\mathbf{\hat X}_n}(\mathbf{\bar X}) = 
\mathbf{\hat X}_n + \mathbf V ~
\sqrt{n} \left(\Psi \left(\frac{\mathbf{
\hat X}_n}{\sqrt n} + \mathbf{\bar X} \right) - \Psi(\mathbf{\bar X}) \right).
\]
By Lemma \ref{lem:X-hat} and the Skorokhod representation theorem
\cite{Bi68}, it follows that 
\begin{equation}\label{eq:srt}
\mathbf{\hat X}_n \to \mathbf{\hat X}~\text{in}~(\sC^K,U)~\text{a.s.}
\end{equation}
as $n \to \infty$. The lemma is proved once we
show that 
\[
\sqrt n \left( \Psi \left(\frac{\mathbf{\hat X}_n}{\sqrt
    n} + \mathbf{\bar X} \right) - \Psi \left(\frac{\mathbf{\hat X}}{\sqrt
    n} + \mathbf{\bar X} \right)\right)  \to 0 ~\text{in}~(\sC^K,U)~\text{a.s. as}~ n \to\infty.
\]
Chen and Whitt \cite{ChWh1993} show that the oblique reflection map
and the reflection regulator are Lipschitz  continuous with respect to
the uniform metric topology. Recall that $\|\cdot\|$ represents the uniform
metric on $\sC^K$ over the interval $[0,T]$. Then,
\begin{align*}
\left\|  \sqrt n \left( \Psi \left(\frac{\mathbf{\hat X}_n}{\sqrt
    n} + \mathbf{\bar X} \right) - \Psi \left(\frac{\mathbf{\hat X}}{\sqrt
    n} + \mathbf{\bar X} \right) \right) \right\|
&\leq K \left\|\sqrt n \frac{\mathbf{\hat X}_n}{\sqrt
    n} + \mathbf{\bar X} - \frac{\mathbf{\hat X}}{\sqrt
    n} - \mathbf{\bar X} \right \|,\\
&= K \left\| \mathbf{\hat X}_n -
    \mathbf{\hat X} \right\|,
\end{align*}
where $K$ is the Lipshitz constant associated with the oblique
reflection map. The conclusion follows as a consequence of~\eqref{eq:srt}.
\EndProof

\Proof [Proof of Lemma~\ref{lem:extra-terms}]
The condition on $x_n, x$ implies that $x_n \stackrel{a.s.}{=} x + (\sqrt n)^{-1} \chi
+ o(\sqrt n)$. Therefore, it follows that
\begin{align*}
    \left \| \sqrt n \left( \Psi(x_n) - \Psi \left( \frac{\chi}{\sqrt
            n} + x\right) \right) \right\| &\stackrel{a.s.}{=}\| \sqrt n \left( \Psi \left( \frac{\chi}{\sqrt
            n} + x + o(1)\right) - \Psi \left( \frac{\chi}{\sqrt
            n} + x\right)\right) \|\\
  &\leq K \sqrt n \| o (1)\|, 
\end{align*}
where the last inequality follows from the Lipshitz continuity of the
oblique reflection map. The final conclusion follows from the fact
that the indeterminate form on the right hand side converges to 0 as
$n \to \infty$.
\EndProof

\noindent \textsc{Remarks:}~We include a short summary of the relevant results
in \cite{MaRa10} that imply that process-level convergence might be
near impossible to prove (in general) in a transitory queueing network. Lemma 2 in
\cite{HoJaWa2012} (an extension of Theorem 3.2 in \cite{MaMa95})
proves the process-level diffusion limit result in the $M_1$ topology
for a single queue. The fact that
the limit process has right- or left-discontinuity points that are `unmatched' by the
pre-limit process necessitates that convergence be proved in the $M_1$
topology as opposed to the more natural $J_1$ topology. On the other
hand, \cite{MaRa10} show that it is not possible to prove a
process-level convergence result even in the $WM_1$ topology (`weak'
$M_1$ topology (see
\cite{Wh01}), due to the fact that the multidimensional limit
process can have discontinuity points that are \textit{both right- and
  left-discontinuous}. For completeness, we state the relevant portion
ofTheorem 1.2
of \cite{MaRa10} that encapsulates the various necessary conditions
for discontinuities in the sample paths of the directional derivative
limit process, $\mathbf
\D_{\mathbf{\hat X}}(\mathbf{\bar X})$. First, given $(z,y)$ as the
solution to the oblique reflection problem for $x \in \sC_0$ define, for each $t \in
[0,\infty)$, 
\begin{eqnarray*}
\sO (t) &:=& \{i \in \{1,\ldots,K\} : z^i(t) > 0\},\\
\sU (t) &:=& \{ i \in \{1,\ldots,K\} : z^i(t) = 0, ~\D y^i(t+) \neq 0,
~\D y^i(t-) \neq 0\}, \\
\sC (t) &:=& \{ 1,\ldots,K \} \backslash [\sO(t) \cup \sU(t)],\\
\sE \sO (t) &:=& \{i \in \sC (t) : \exists \delta > 0 \text{ such that } z^i(s)
> 0 ~\forall s \in (t-\delta,t)\},\\
\sS \sU (t) &:=& \{ i \in \sC (t) : \D z^i(t-) = 0, ~\D z^i(t+) \neq 0.
\}
\end{eqnarray*}
When $x = \mathbf{\bar X}$, $\sO (t)$ is the set of nodes in the
network that are \emph{overloaded} at time $t$, $\sU (t)$ is the set of
underloaded nodes, $\sC (t)$ the set of critically loaded nodes, $\sE
\sO (t)$ is the set of critically loaded queues that are at the end of
overloading and $\sS \sU (t)$ is the set of critically loaded nodes
that are at the start of under-loading. Note that the definitions of
overloading, under-loading and critical loading conform to the standard
notions for $G/G/1$ queues, as noted in \cite{HoJaWa2012}. Next, we
also require the notion of critical and sub-critical chains, as in
Definition 1.5 of \cite{MaRa10}:

\begin{definition}[Def. 1.5~\cite{MaRa10}]\label{def:mara}
  Given a $K \times K$ routing matrix $\mathbf P$ and the oblique reflection map
  $\Psi$ and $x \in \sC^K$ so that $y = \Psi(x)$. Then a sequence
  $j_0, j_1, \ldots, j_m$ with $j_k \in \{1,\ldots,K\}$ for $k =
  0,1,\ldots,m$ that satisfies $P_{j_{k-1},j_k} > 0$ for $k =
  0,1,\ldots,m$ is said to be a chain. The chain is said to be a cycle
  if there exist distinct $k_1, k_2 \in \{0,\ldots,m\}$ such that
  $j_{k_1} = j_{k_2}$, the chain is said to precede $i$ if $j_0 = i$
  and is said to be empty at $t$ if $y_{j_k}(t) = 0$ for every
  $k=1,\ldots,m$. For $i=1,\ldots,K$ and $t \in [0,\infty)$, we
  consider the following two types of chains:
  \begin{enumerate}
  \item An empty chain preceding $i$ is said to be critical at time
    $t$ if it is either cyclic or $j_m$ is at the end of overloading
    at $t$.
    \item An empty chain preceding $i$ is said to be sub-critical at time
      $t$ if it is either cyclic or $j_m$ is at the start of
      overloading at $t$.
  \end{enumerate}
\end{definition}

Theorem 1.2
of \cite{MaRa10} gives necessary conditions so that, in general, the sample paths of the
directional derivative can have both a right \emph{and} left
discontinuity at $t \in [0,\infty)$. Simply put, the structure of the
routing matrix $\mathbf P$ determines whether we see such a point. 
\begin{proposition}[Thm. 1.2~\cite{MaRa10}] \label{prop:mara10}
  Under the conditions of Definition~\ref{def:mara} and given a
  process $\chi \in \sC^k$, if the directional
  derivative $\D_{\chi}(x)$ has both a right and a left discontinuity
  at $t \in [0,\infty)$ then one of the following conditions must hold
  at time $t$:
  \begin{enumerate}
  \item[a)] $i$ is at the end of overloading, and a sub-critical chain
    precedes $i$, in which case
    \[
    \D_{\chi}(x)^i(t-) < \D_{\chi}(x)^i(t)^i = 0 < \D_{\chi}(x)^i(t+),
    \]
    \item[b)] $i$ is at the start of under-loading and a critical chain
      precedes $i$, in which case 
      \[
      \D_{\chi}(x)^i(t-) > \D_{\chi}(x)^i(t) > \D_{\chi}(x)^i(t+) = 0,
      \]
      \item[c)] $i$ is not underloaded and there exist both critical
        and sub-critical chains preceding $i$; if, in addition, $i$ is
        overloaded then the discontinuity is a separated discontinuity
        of the form
        \[
        \D_{\chi}(x)^i(t) < \min\{\D_{\chi}(x)^i(t-), \D_{\chi}(x)^i(t+)\}.
        \]
  \end{enumerate}
\end{proposition}
Note that the sample paths of $\mathbf
\D_{\mathbf{\hat X}}(\mathbf{\bar X})$ lie in $\sD^K_{\lim}$ and
establishing $M_1$ convergence in this space is non-trivial. Recall that
the standard description of $M_1$ convergence is through the graphs of
the functions - which can be described via linear interpolations in
$\sD$ and $\sD^K_{l,r}$. However, in $\sD^K_{\lim}$ no such simple
description exists (see Chapter 12 of \cite{Wh01} and Chapter 6, 8 of
\cite{Wh01b} for further details on these issues). 

Given the inherent difficulty in establishing a general process-level result, we
first focus on a two queue tandem network, where the arrival time
distribution is uniform on the interval $[-T_0,T]$ and $T_0, T > 0$ where the difficulties will become
apparent.

\begin{theorem} \label{thm:tandem-uni-diffusion}
Consider a tandem queueing network with 
\(
\mathbf P = \begin{pmatrix} 
0 & 0 \\
1 & 0
\end{pmatrix},
\)
and $\mathbf R = \mathbf I - \mathbf P^T = $
\(
\begin{pmatrix}
1 & 0\\
-1 & 1
\end{pmatrix}.
\) 
Assume that $\mathbf F = F_1$ is uniform over $[-T_0,T]$, and service rate at node 1
is $\m_1$ and at node 2 $\m_2$. Then, 
\(
\mathbf{\hat Q}_n \Rightarrow \mathbf{\hat Q} := \mathbf \D_{\mathbf{\hat X}}(\mathbf
{\bar X}) 
\)
in $(\sD^2_{l,r},SM_1)$ as $n \to \infty$, where
$\mathbf{\hat X} = (\hat X_1, \hat X_2)$ with $\hat X_1 = W_1^0 \circ F_1
- W_1 \circ M_k$, $\hat X_2 = W_1 \circ M_k - W_2 \circ M_2$ and
$M_k(\cdot) = \int_0^\cdot \mu_k(s) ds$ for $k \in \{1,2\}$, $\mathbf{\bar X} = ((F_1 - \m_1 e), (\m_1 -
\m_2) e)^T$ and $e :\bbR \to \bbR$ is the identity map.
\end{theorem}
\Proof
Recall that $F(t) =
\frac{t+T_0}{T+T_0}$ for all $t \in [-T_0,T]$. We consider three subcases and establish the weak convergence result
for each of them separately.

\noindent (i) Let $\m_1 < \m_2$. Then,
\begin{equation} \label{node1-fluid}
\bar Q_1(t) = \begin{cases}
(F(t) - \m_1 t \mathbf 1_{\{t \geq 0\}}) & \forall t \in [-T_0,
\tau_1),\\
0 & \forall t \in [\tau_1,\infty),
\end{cases}
\end{equation}
and 
\(
\bar Q_2(t) = 0 ~\forall t \geq 0,
\)
where $\tau_1 := \inf\{t > 0 | F(t) = \m_1 t\}$. These follow as a
consequence of Corollary \ref{cor:tandem-fluid}, and noting that $\mathbf{\bar
X} = (F(t) - \m_1 e, (\m_1 - \m_2) e)$. Thus, we have
\begin{eqnarray}
\label{node1-nabla}
\nabla_t^1 & := & \begin{cases}
\{-T_0\} & \forall t \in [0,\tau_1),\\
\{-T_0, \tau_1\} & t = \tau_1,\\
\{t\} & \forall t > \tau_1, \text{ and}
\end{cases} \\
\label{node2-nabla}
\nabla_t^2 &:=& \{t\} ~\forall t \in [0,\infty).
\end{eqnarray}
Thus, node 1 is in $\sO (t)$ for all $t \in [-T_0,\tau_1)$, $\sC (t)$
for $t = \tau_1$ and in $\sU (t)$ for $t > \tau_1$, and node 2 is 
in $\sU (t)$ for all $t$. 

The limit process $\mathbf{\hat Q}$ has a
discontinuity only in the first component at $\hat Q_1(\tau_1) = \hat
X^1(\tau_1) + \max\{0,-\hat X^1(\tau_1)\}$. Note that $\hat Q_1(\tau_1 -) =
\hat X^1(\tau_1)$ and $\hat Q_1(\tau_1+) = 0$, implying that $\hat Q_1$
has either a right or left discontinuity at $\tau_1$. If $\hat
X^1(\tau_1) \geq  0$ then $\hat Q_1(\tau_1) = \hat X^1(\tau_1) = \hat Q_1(\tau_1
-) > \hat Q_1(\tau_1+) = 0$ and has a right discontinuity. Else, if $\hat X^1(\tau_1) < 0$ then $\hat
Q_1(\tau_1) = 0 = \hat Q_1(\tau_1+) > \hat Q_1(\tau_1 -)$ and has a left
discontinuity. Thus, the limit process $\mathbf{\hat Q}$ has sample
paths in $\sD^2_{l,r}$. The proof of convergence for $\mathbf{\hat Q}_n
= (\hat Q_{n,1}, \hat Q_{n,2})$ in this case is
simple. First, Theorem 2 of \cite{HoJaWa2012} shows that $\hat Q_{n,1}
\Rightarrow \hat Q_1 := \hat X^1 + \sup_{s \in \nabla_{\cdot}}^1
(-\hat X(s))$ in $(\sD_{l,r}, M_1)$ as $n \to \infty$, and
$\hat Q_{n,2} \Rightarrow 0$ in $(\sD_{l,r}, M_1)$. Recall that $Disc(\hat
Q_1)$ and $Disc(\hat Q_2)$ are the (respective) sets of discontinuity
point, and it is obvious that $Disc(\hat Q_1) \cap Disc(\hat Q_2) =
\phi$. Therefore, by \cite[Corollary 6.7.]{Wh01b}, 
\(
\hat Q_{n,1} + \hat Q_{n,2} \Rightarrow \hat Q_1
\)
in $(\sD_{l,r}(\bbR),M_1)$ as $n \to \infty$. Consequent to \cite[Theorem 6.7.2]{Wh01b}, it follows that
$\mathbf{\hat Q}_n \Rightarrow \mathbf{\hat Q} := (\hat Q_1, 0)^T$ in
$(\sD^2_{l,r}, SM_1)$ as $n \to \infty$. 

\noindent (ii) Let $\m_1 > \m_2$. Then, $\bar Q_1$ and $\nabla_t^1$
follow \eqref{node1-fluid} and \eqref{node1-nabla} (resp.). $\bar Q_2$
on the other hand, is more complex now:
\[
\bar Q_2(t) = \begin{cases}
(\m_1 - m_2) t & \forall t \in[0,\tau_1],\\
(F_1(t) - \m_2 t) &
\forall t \in [\tau_1,\tau_2],\\
0 & \forall t > \tau_2,
\end{cases}
\]
where $\tau_2 := \inf \{t > \tau_1 : F_1(t) = \m_2 t\}$ (note that $\tau_2 >
\tau_2$ since $\m_1 > \m_2$). It follows that
\[
\nabla_t^2 = \begin{cases}
\{0\} & \forall t \in [0,\tau_2),\\
\{0,\tau_2\} & t = \tau_2,\\
\{t\} & \forall t > \tau_2.
\end{cases}
\]
It follows that node $2$ is in $\sO(t)$ for all $t \in [0,\tau_2)$,
$\sC(t)$ at $t = \tau_2$ and $\sU(t)$ for $t > \tau_2$.

The diffusion limit $\mathbf{\hat Q} := (\hat Q_1, \hat Q_2)$ has discontinuities in both
components. For node 1, if $\hat X^1(\tau_1) \geq 0$ then $\hat
Q_1(\tau_1)$ has a right discontinuity, while $\hat X^1(\tau_1) < 0$ then
$\hat Q_1(\tau_1)$ has a left discontinuity. Similarly, if $\hat
X^2(\tau_2) \geq 0$ then $\hat Q_2(\tau_2)$ has a right discontinuity, and
if $\hat X^2(\tau_2) < 0$ it has a left discontinuity. It follows that
$\mathbf{\hat Q}$ has sample paths in $\sD^2_{l,r}$. Furthermore, it is
clear that $Disc(\hat Q_1) \cap Disc(\hat Q_2) = \phi$. Therefore, the
weak convergence result follows by the same reasoning as in part (i).

\noindent (iii) Assume $\m_1 = \m_2$. Once again, $\hat Q_1$ and
$\nabla_t^1$ follow \eqref{node1-fluid} and \eqref{node1-nabla}
(resp.). On the other hand, for node 2 $\hat Q_2 = 0$, but unlike case
(i), the queue is empty but the server operates at full capacity till
$\tau_1$, and then enters underload. Thus,
\[
\nabla_t^2 = \begin{cases}
[0,t] & \forall t \in [0,\tau_1],\\
\{t\} & \forall t > \tau_1.
\end{cases}
\]
It is clear that node 2 switches from $\sC(t)$ in $[0,\tau_1]$ to
$\sU(t)$ for $t > \tau_1$. Furthermore, at $\tau_1$ itself, the node is in
$\sS \sU(t)$ (the regulator is flat to the left of $\tau_1$ and
increasing to the right).

The diffusion limit, once again, has discontinuities in both
components. However, it is clear that $Disc(\hat Q_1) = Disc(\hat
Q_2) = \{\tau_1\}$. For any $\sT > -T_0$, it is straightforward to see
that $(\hat Q_1(t) - \hat Q_1(t-)) (\hat Q_2(t) - \hat Q_2(t-)) \geq
0$ for all $-T_0 \leq t \leq \sT$: clearly, for any $t < \tau_1$, $\hat Q_i$, $i = 1,2$ are both
continuous. On the other hand, at $\tau_1$, $\hat Q_1(\tau_1) \geq \hat
Q_{1}(\tau_1-)$ and $\hat Q_2(\tau_1) = \hat Q_2(\tau_1-)$. Finally, for any
$t > \tau_1$, $\hat Q_1(\tau_1) = \hat Q_{1}(\tau_1-)$ and $\hat Q_2(\tau_1) =
\hat Q_2(\tau_1-)$. Now, by Theorem 6.7.3 of \cite{Wh01b}, it follows
that $\hat Q_{n,1} + \hat Q_{n,2} \Rightarrow \hat Q_1 + \hat Q_2$ in
$(\sD_{l,r}(\bbR),M_1)$ as $n \to \infty$. Then, by Theorem 6.7.2 of
\cite{Wh01b}, $\mathbf{\hat Q}_n \Rightarrow \mathbf{
\hat Q}$ in $(\sD^2_{l,r},SM_1)$ as $n \to \infty$. This concludes the proof.
\EndProof

Theorem \ref{thm:tandem-uni-diffusion} shows that in the case of a
tandem network, with uniform arrival time distribution, the weak
convergence result can be established in the space $\sD^2_{l,r}$ and in
the $SM_1$ topology. In fact this result is true, if $F_1$ is unimodal
such that node 1 is overloaded in the initial phase (i.e., in the
interval $[-T_0,\tau_1)$, with $T_0 \geq 0$ now). We capture this fact
in the following corollary. Without loss of generality, we will assume
that $T_0 = 0$.

\begin{corollary} \label{lem:tandem-arbit-diffusion}
Let $F_1$ be a unimodal distribution function with finite support
$[0,T]$, and consider a tandem queue as defined in Theorem
\ref{thm:tandem-uni-diffusion}. Then, $\mathbf{\hat Q}_n \Rightarrow
\mathbf{\hat Q} := \mathbf{\D}_{\mathbf{\hat X}}(\mathbf{\bar X})$ in
$(\sD^2_{l,r},SM_1)$ as $n \to\infty$, where 
\[
\mathbf{\hat X} := \left(
  W_1^0
\circ F_1 - \s_1 \m_1^{3/2} W_1, (\s_1 \m_1^{3/2} W_1 - \s_2 \m_2^{3/2}
W_2) \right)^T,
\] 
$\mathbf{\bar X} = (F_1 - \m_1 e, (\m_1 - \m_2)e)^T$
and $e : \bbR \to \bbR$ is the identity map. 
\end{corollary}

The proof follows that of Theorem \ref{thm:tandem-uni-diffusion} and
is omitted. Note that the compact support assumption is required, due
to the fact that we prove weak convergence over compact intervals of
time (see Section 7.2 of \cite{HoJaWa2012} for a discussion on this
point).

\section{High-intensity Analysis of Tandem Networks}
We illustrate the utility of the afore-developed approximations in
bottleneck analysis of transitory tandem networks. Almost all of the analysis
in the literature has focused on the characterization and detection of
bottlenecks in stationary queueing networks. Of particular relevance
to our results in this paper is the heavy-traffic bottleneck
phenomenon identified in \cite{SuWh1990,Wh01}. To recall, the \textit{heavy-traffic
bottleneck phenomenon} corresponds to the state space collapse that
is observed when the traffic intensity at a single queue
approaches 1, while the traffic intensity at other queues remains
below 1. In this case, the well known heavy-traffic approximations in
\cite{IgWh70b,Re1984,ChMa1991} imply that the network workload
process will collapse to a single dimensional process determined by the
bottleneck node. In other words, the non-bottleneck nodes behave like
`switches' where the service time is effectively zero.~%  In
% \cite{SuWh1990} the authors describe a simulation experiment that
% illustrates that the heavy-traffic phenomenon can be significant even
% at reasonably moderate traffic intensities. The authors consider a
% series queue network fed by a unit-rate renewal traffic process, and
% exponential service at each node in the network (that is also
% independent of the service at the other nodes). Under the assumption
% that the final node in the network alone approaches heavy-traffic, it
% is shown that the steady-state waiting time distribution at the
% bottlneck node is asymptotically (as the traffic intensity approaches 1) equal to
% the steady-state waiting time distribution if the departure process
% from the previous node were replaced by the renewal traffic process
% feeding into the network. 
In general, exact bottleneck analysis is very difficult (if
not impossible), and several approximations have been
proposed in the literature, particularly the parametric-decomposition
approach \cite{Wh1983,BuSh1992}, the stationary-interval
method \cite{Wh1984}, and Reiman's individual (IBD) and sequential bottleneck
decomposition (SBD) algorithms \cite{Re1990}. % As noted in \cite{Re1990}, IBD is
% equivalent to the asymptotic method introduced in \cite{Wh1984}. The
% key idea behind the IBD is to treat each node in the network as the
% unique bottleneck, and all other nodes as instantaneous
% switches. 
Nonetheless, the natural metric to use to study bottlenecks would be
the waiting time at each node. The fluid and diffusion workload
approximations can be established as a corollary to Theorem
\ref{thm:queue-fluid} and Theorem~\ref{thm:pt-wise}, assuming that the service process is stationary:

\begin{corollary} [Workload Approximation] \label{prop:workload}
  Recall that $\mathbf M$ is a diagonal matrix defined as 
\[
\mathbf M := \text{diag}(1/\mu_1,\ldots, 1/\mu_1,1/\mu_K). 
\]
Then the fluid workload process $\mathbf {\bar  Z} =
  \mathbf {M \bar Q}$, and for each $t \in [0,\infty)$
  the diffusion workload process is $\mathbf{\hat Z}(t) = \mathbf{M \hat Q}(t)$.
\end{corollary}
The proof of this corollary follows by analogous arguments to
\cite[Proposition 4]{HoJaWa2012}, and we omit it. 

Bottleneck analysis, however, has largely
been ignored in transitory networks
in particular. The key difference (and
difficulty) in the transitory setting is that, for general arrival
epoch distributions $F$ it is possible that the number of bottleneck queues can
change with time. The
situation is considerably simpler when $F$ is uniform, however, and we focus on
this case first to illustrate the main ideas. We commence with
a definition of a bottleneck queue in a transitory network, in the
large population limit.

\begin{definition} [Transitory Bottleneck Queue]
  Queue $k \in \sK$ in the transitory queueing network is a bottleneck at time $t$ if and
  only if the diffusion workload process satisfies $|\hat Z_k(t)| > 0$. 
\end{definition}

Note that we choose to use a sample path definition of the bottleneck
node owing to the fact that the temporal stochastic variations can produce
differing numbers of bottlenecks, even compared with the average/fluid
variation. This definition is natural to consider in job-shop type
production systems and complements the definitions
in~\cite[P. 23]{LaBu1994} that classifies bottlenecks in terms of short,
intermediate and long time horizons. 

\begin{example} [Tandem network with uniform traffic]~\label{ex:uniform}
Consider a series network of $K$ queues. Let the service rate
at queues $1$ through $K-1$ be $\mu_1$ and $\mu_{K}$ at queue
$K$. Without loss of generality we assume that $\mu_K < 1 \leq
\mu_1$. Assume that the traffic arrival epochs are
randomly scattered per a uniform distribution function, over the
interval $[0,1]$.  Then, in the fluid population acceleration limit as
observed in Theorem ~\ref{thm:queue-fluid}, it
can be observed that each of the queues $1,\ldots,K-1$ behave like
instantaneous switches and $O(n)$ fluid accumulates at the final
queue. Extending the analysis in Corollary~\ref{cor:tandem-fluid} to a
$K$-node tandem network it is straightforward to compute that $\bar {\mathbf{X}} =
(\bar X_1,\ldots, \bar X_{K-1},\bar X_K)$, where $\bar X_1(t) = F_1(t) - \mu_1 t \leq 0$ and
$\bar X_k(t) = 0$ for all $k = 2,\ldots, K-1$, and
$\bar X_K(t) = (\mu_{1} - \mu_K) t > 0$. Since the routing matrix is
\[
\mathbf{P} = 
\begin{pmatrix}
  0 & 0  & \ldots & 0 & 0\\
  1 & 0 & \ldots & 0 & 0\\ 
  0 & 1 & \ldots & 0 & 0\\ 
  & & \vdots &  &\\
  0 & 0 & \ldots & 1 & 0\\  
\end{pmatrix}
\]
a simple (if tedious) calculation shows that 
\[
\bar{\mathbf{ Q}}(t) =
\begin{cases}
  \left( 0,\ldots,0, (\mu_1- \mu_K)t \right) & t \in [0,1/\mu_K],\\
  \left(0,\ldots, 0\right) & t > 1/\mu_K.
\end{cases}
\]
Now, It follows that,
in the case of the tandem queueing network under consideration
\[
\bar{\mathbf{Z}}(t) =
\begin{cases}
  (0, \ldots, 0, (\mu_1\mu_K^{-1} - 1)t) & t \in [0,1/\mu_K],\\
  (0,\ldots, 0) & t > 1/\mu_K.
\end{cases}
\]
Thus, in the fluid limit, we find
that the tandem queueing network ``collapses'' to a single queue
in the fluid limit (this is an example of a \textit{state
  space collapse} as defined in \cite{Re1984}), and the sojourn time
through the network, in the fluid scale and large population limit, is
determined entirely by the delay at node $K$. 

On the other hand, as the diffusion limit
in Theorem \ref{thm:tandem-uni-diffusion} shows, there is non-zero
variability in the queue length at each node in the network. Indeed,
Theorem \ref{thm:tandem-uni-diffusion} and
Corollary~\ref{prop:workload} imply that the diffusion limit of the
workload vector in a tandem network is $\mathbf {\hat Z}
= \mathbf{M}\mathbf{\D}_{\hat{\mathbf{X}}}(\bar{\mathbf{X}})$, where
\[
\begin{split}
\hat{\mathbf{X}}(t) = \bigg((W_1^0(t)- \sigma \mu_1^{3/2} W_1(t)), (\sigma_1
  &\mu_1^{3/2} W_1(t) - \sigma_1 \mu_1^{3/2} W_2(t), \ldots,\\
  &\sigma_1\mu_1^{3/2}W_{K-1}(t) - \sigma_K \mu_K^{3/2} W_K(t))
\bigg).
\end{split}
\]
 Now, if
$\mu_1 > 1$, then $\hat Z_k \stackrel{D}{=} 0$ for $k = 1,\ldots,K-1$ and
$\hat Z_K(t) \stackrel{D}{=} \mu_K^{-1} (\hat X_K(t) + \sup_{0 \leq s \leq
  t} (-\hat X_K(s)))$ with $\hat X_K = \sigma_1\mu_1^{3/2}W_{K-1} - \sigma_K
\mu_K^{3/2} W_K$. That is, in the population
acceleration scaling, the distribution of the sojourn time through the
network is asymptotically equal to the delay distribution of the last queue.

On the other hand,
if $\mu_1 = 1$, then $\hat Z_1 = \mu_1^{-1} (\hat X_1(t) + \sup_{0 \leq s
  \leq t} (-\hat X_1(s)))$ with $\hat X_1 = W_1^0 - \sigma \mu_1^{3/2}
W_1$, $Z_k \stackrel{D}{=} 0$ for $k=2,\ldots, K-1$
and 
\[
\hat Z_K = \begin{cases}
  \mu_K^{-1}(\sigma_1 \mu_1^{3/2} W_{K-1} - \sigma_K \mu_K^{3/2} W_K)
  ~&\forall t \in [0,1]\\
  \mu_K^{-1} (- \sigma_K \mu_K^{3/2} W_K) ~&\forall t \in (1,1/\mu_K]\\
 0 ~&\forall t > 1/\mu_K.
\end{cases}
\]

This indicates that there are two bottlenecks at queues $1$ and
$K$. Thus, there is a state space collapse to a two-dimensional
vector $\mathbf {\hat Z} = (Z_1, Z_K) $, and the sojourn time through
the network is
asymptotically equal in distribution to the sum of the delays in these
two queues. 
\end{example}

\begin{example}[Tandem network with unimodally traffic] \label{ex:unimodal}
Now, suppose $F_1$ is not uniform, but unimodal with support on
$[0,1]$.  For
simplicity, we assume that the distribution function is symmetric around
$\tau:= \text{argmax}\{F_1'(t) : t \in [0,1]\}$, where $F_1'$ represent the
first derivative of the arrival epoch distribution (assuming it is
well defined), and that the service
rates are the same as in Example~\ref{ex:uniform}. The uni-modality of the
arrival epoch distribution implies that up to time $\tau$ the distribution function is convex
increasing, while after $\tau$ it is concave decreasing. A simple
example of such a distribution function would be,
\begin{align*}
  F_1'(t) = \begin{cases}
    4 t ~ & t \in [0,1/2]\\
    4(1-t) ~& t \in (1/2,1].
    \end{cases}
\end{align*}
In this case, $\tau = 1/2$ and $F_1'(\tau) = 2$. 

% If the service rates satisfy $\mu_K < F_1'(\tau) \leq \mu_1$, then the
% fluid and diffusion queue length processes are exactly as described in
% Example~\ref{ex:uniform}. 

We first focus on the case where the service rates satisfy
$\mu_K < \mu_1 < F_1'(\tau)$. Observe that 
\[
\bar X'_1(t) = F'_1(t) - \mu_1 \begin{cases}
  \leq 0 ~&\forall t \in [0,\tau_1)\\
  > 0 ~&\forall t \in [\tau_1,\tau_2),\\
  \leq 0 ~& \forall t \in [\tau_2,1],
\end{cases}
\]
where $\tau_1 := \inf\{t > 0 : F'(t) = \mu_1\}$ and $\tau_2 := \inf\{t
> \tau_1 : F'(t) = \mu_1\}$; that is, these are the two points in time
where the derivative of the arrival epoch distribution equals the
service rate in queue 1. We also have $\bar X_k(t) =
0$ for all $k = 2,\ldots,K-1$ and $\bar X_K(t) = (\mu_1 - \mu_K)t >
0$, for all $t \in [0,1]$. Consider the fluid queue length at node 1 $\bar Q_1(t) = \bar
X_1(t) + \sup_{0\leq s \leq t}[-\bar X_1(s)]^+$, and observe that 
\begin{align} \label{eq:node1-f}
  \bar Q_1(t) =
  \begin{cases}
    0 ~& t \in [0,\tau_1)\\
    \bar X_1(t) - \bar X_1(\tau_1)~& t \in [\tau_1, \tau_2)\\
    0 ~& t \geq \tau_2.
  \end{cases}
\end{align}
% Observe that, unlike in Example~\ref{ex:uniform}, queue 1 is a
% bottleneck in the interval $[\tau_1,\tau_2)$. 
Following the arguments in
Example~\ref{ex:uniform} it can be shown that the fluid
queue length in the downstream nodes satisfies $\bar Q_k(t) = 0$ for
all $k=2,\ldots,K-1$ for all $t > 0$. Similarly, in the terminal node
\begin{align*}
\bar Q_K(t) = 
  \begin{cases}
    0 ~& t \in [0,\tau_1')\\
    (F_1(t) - F_1(\tau_1')) - \m_K(t - \tau_1') ~& t \in
    [\tau_1',\tau_1)\\
    (F_1(\tau_1) - F_1(\tau_1')) - \mu_1(\tau_1 - \tau_1') + (\mu_1-\mu_K) (t
    - \tau_1') ~& t \in [\tau_1,\tau_2)\\
    (F_1(t) - F_1(\tau_1')) - \mu_K(t - \tau_1') - (F_1(\tau_1) - \mu
    \tau_1) ~& t \in [\tau_2,\tau_2')\\
    0 ~& t \geq \tau_2',
  \end{cases}
\end{align*}
where $\tau_1' := \inf \{t > 0 : F_1'(t) \geq \mu_K\}$ and $\tau_2' :=
\sup \{t > \tau_1' : F_1'(t) \geq \mu_K\}$. This follows from the fact
that
\[
F_1'(t) - \mu_K \begin{cases}
  \leq 0 ~&\forall t \in [0,\tau_1')\\
  > 0 ~&\forall t \in [\tau_1',\tau_2'),\\
  \leq 0 ~& \forall t \in [\tau_2',1],
\end{cases}
\]
% To summarize, the fluid queue length process satisfies,
% \begin{align*}
%   \bar{\mathbf Q}(t) = 
%   \begin{cases}
%     \left(0, 0, \ldots,0,(\mu_1 - \mu_{K})t\right) ~& t \in
%     [0,\tau_1)\\
%     \left((\bar X_1(t) - \bar X_1(\tau_1), 0, \ldots, 0, (\mu_1 -
%       \mu_{K})t) \right) ~& t \in [\tau_1, \tau_2)\\
%     \left(0,0,\ldots (\mu_1 - \mu_K)t\right) ~& t \in [\tau_2,1/\mu_K]\\
%     \left(0,\ldots,0 \right) ~& t > 1/\mu_K.
%   \end{cases}
% \end{align*}
% Now, as a consequence of Corollary~\ref{prop:workload}, the fluid workload
% process is
% \[
%  \bar{\mathbf Z}(t) = \begin{cases}
%     \left(0, 0, \ldots,0, 0\right) ~& t \in
%     [0,\tau_1')\\
%     \left(0, 0, \ldots,0, \right) ~& t \in
%     [0,\tau_1')\\
%     \left(\mu_1^{-1}(\bar X_1(t) - \bar X_1(\tau_1), 0, \ldots, 0,
%       (\mu_1 \mu_{K}^{-1} - 1)t) \right) ~& t \in [\tau_1, \tau_2)\\
%     \left(0,0,\ldots (\mu_1 \mu_K^{-1} - 1)t\right) ~& t \in [\tau_2,1/\mu_K]\\
%     \left(0,\ldots,0 \right) ~& t > 1/\mu_K.
%   \end{cases}
% \]
In contrast to Example~\ref{ex:uniform} the state-space
collapse is not straightforward here. The fluid tandem queueing network
switches between collapsing to a single queue network in time
intervals $[\tau_1',\tau_1)$ and $[\tau_2, \tau_2')$ and a two queue network
in the interval $[\tau_1,\tau_2)$. Thus, the state-space collapse \textit{itself} exhibits non-stationary behavior.

Next, considering the diffusion limit, extending 
Corollary~\ref{lem:tandem-arbit-diffusion} to a $K$-node network we have
\[
\begin{split}
  \hat{\mathbf X} = \bigg( (W_1^0 \circ F_1 - \sigma \mu_1^{3/2} W_1)&
    \sigma_1 \mu_1^{3/2} (W_1 - W_2),\\ &\ldots, (\sigma_1 \mu_1^{3/2}
    W_{K-1} - \sigma_K \mu_K^{3/2} W_K) \bigg),
\end{split}
\]
and from Corollary~\ref{prop:workload} the workload diffusion limit
process is $\hat{\mathbf Z} = \mathbf M \mathbf{\D}_{\hat{\mathbf
    X}}(\bar{\mathbf X})$. Note that $\hat X_k := \sigma_1 \mu_1^{3/2} (W_{k-1} - W_{k})
\stackrel{D}{=} \sigma_1 \mu_1^{3/2} W_k^*$ for $k = 2,\ldots,K-1$
where $W_k^*$ are independent but identically distributed Brownian
motion processes. 

% Paralleling Example~\ref{ex:uniform}, when $\mu_1 \geq F_1'(\tau)$
% $\hat Z_k \stackrel{D}{=} 0$ for $k = 1, \ldots, K-1$ and $\hat
% Z_K(t) = \mu_K^{-1} \left( \hat X_K(t) + \sup_{0 \leq s \leq t}
%   (-\hat X_K(s)) \right)$. That is, in the diffusive scaling, the
% distribution of the sojourn time is determined entirely by the delays
% experienced at queue $K$.
Following the fluid
limit discussion above, the diffusion limit workload
process at queue 1 satisfies
\[
\hat Z_1(t) = \begin{cases}
  0 ~& t \in [0,\tau_1)\\
  \mu_1^{-1} \left( \hat X_1(t) - \hat X_1(\tau_1) \right) ~& t \in
  [\tau_1, \tau_2)\\
  0 ~& t \geq \tau_2.
\end{cases}
\]
On the other hand, following Corollary~\ref{prop:workload} and using
the description of the directional derivative regulator map in
Lemma~\ref{lem:dir-der-ORT}, the processes $\hat Z_k$ for $k = 2,
\ldots, K-1$ can be shown to satisfy
\begin{align*}
  \hat Z_k(t) =
  \begin{cases}
    0 ~& t \in [0,\tau_1)\\
    \mu_1^{-1} \left(\hat X_k(t) + \sup_{\tau_1 \leq s \leq t} [-\hat
      X(s)]^+ \right) ~& t \in [\tau_1, \tau_2)\\
    0 ~& t \geq t \geq \tau_2.
  \end{cases}
\end{align*}
Observe that jobs flowing through queues $2,\ldots,K-1$ will
experience non-zero delays in the interval $[\tau_1, \tau_2)$,
determined by the reflected Brownian motion process. The
reason why this happens is manifest: the departure rate from queue 1
reaches its maximum value ($\mu_1$) in this interval, so that the
downstream queues $2,\ldots, K-1$ become critically loaded in this
interval. Note that the jobs experience (effectively) zero service delay
in the latter queues, and they are instantaneously switched through to
downstream nodes. Thus, the ``surge'' period $[\tau_1,\tau_2)$ is the
same in all of these nodes. Finally, in queue $K$ the diffusion limit
workload process satisfies
\begin{align}\label{eq:nodeK-d}
  \hat Z_K(t) =
  \begin{cases}
    0~& t \in [0,\tau_1')\\
    \mu_K^{-1} \hat X_K(t) ~& t \in [\tau_1',\tau_2')\\
    0 ~& t \geq \tau_2'.
  \end{cases}
\end{align}
Unlike
Example~\ref{ex:uniform}, jobs experience delays in the last queue
only in the interval $[\tau_1',\tau_2')$. Note that this interval is
includes the interval $[\tau_1, \tau_2]$, due to the assumption that
$F_1$ is unimodal. 

Now, consider the alternative case where the serve rates satisfy $\m_K <
F_1'(\tau) \leq \m_1$. In the fluid approximation, nodes
$k=1,\ldots,K-1$ are always 'underloaded' and thus, $\bar Q_k = 0$. On
the other hand, the queue length at node $K$ satisfies
\[
 \bar Q_K(t) =
  \begin{cases}
    0 ~& t \in [0,\tau_1')\\
    \bar X_1(t) - \bar X_1(\tau_1')~& t \in [\tau_1', \tau_2')\\
    0 ~& t \geq \tau_2',
  \end{cases}
\]
where, as before, $\tau_1' := \inf \{t > 0 : F_1'(t) \geq \mu_K\}$ and $\tau_2' :=
\sup \{t > \tau_1' : F_1'(t) \geq \mu_K\}$. Consequently, it follows
that the diffusion approximation satisfies $\hat Z_k \stackrel{D}{=}
0$ for $k = 1,\ldots,K-1$ and $\hat Z_K(t)$
follows~\eqref{eq:nodeK-d}. That is, the only bottleneck node in the
network manifests at node $K$ in the interval $[\tau_1',\tau_2')$.
\end{example}

As a final observation, note that the diffusion limit processes all exhibit
discontinuities - at time $\tau_2$ for queues $k = 1, \ldots, K-1$ and
at $\tau_2'$ for queue $K$. This parallels similar observations for
the single server transitory queue in~\cite{HoJaWa2012}. Again unlike Example~\ref{ex:uniform}, the bottlenecks in the network change
over time: there are no bottlenecks in $[0,\tau_1')$, one bottleneck
(queue $K$) in $[\tau_1',\tau_1)$, all the queues are bottlenecks in
$[\tau_1,\tau_2)$, one bottleneck (queue $K$) in $[\tau_2,\tau_2')$
and zero bottlenecks in $[\tau_2',\infty)$. Thus, the state space
``collapse'' of the diffusion appoximation is much more complicated.

\section{Concluding Statements} \label{sec:conclusions}
In this paper we developed asymptotic `population acceleration' approximations of
the queue length and (implicitly) the workload processes in a network
of transitory queues. These results complement and add to the body of
research studying single class generalized Jackson networks. In particular, our fluid limit results
accomodate rather general traffic and service models. On the other
hand, we can only establish point-wise diffusion approximations in the
most general case, owing to the difficulties in the existence of the
so-called directional derivative oblique reflection map. Nonetheless,
we establish functional central limit theorems in the special case of
a tandem network and we also present direct consequences of these
developments on bottleneck analysis. 

There are several directions in which this research will be expanded
in the future. The extension of these results to general polling
queueing networks will be interesting, exploiting some recently
observed connections between acceleration scalings and polling
networks in \cite{Ve2015}. Second, the arrival
counts in non-overlapping intervals under the $\Delta_{(i)}$ traffic
model have strong negative association. How soon will this correlation
be `forgotten' as traffic passes through multiple stages of service?
This requires a study of the possible sample paths of the workload
process. We believe this question has deep connections with directed
percolation models; this is not a novel observation: \cite{GlWh1991} identify this
connection when there are no traffic dynamics. In on-going work we are
working towards extending their analysis to transitory networks. A further interesting question is how the last passage
percolation time scales with the population size in a non-stationary
setting (as opposed to the classical setting where the percolation
model is only studied in the stationary setting). The connection
between percolation time and the sojourn time through the network
affords yet another bottleneck/performance analysis measure in
networks of transitory queues that will be highly relevant in the
context of manufacturing lines. We will consider these questions in
future papers.

\subsection{Appendix subsection}
\noindent\textbf{Proof of Theorem~\ref{thm:arrival-limits}}
We prove lemma's for each of the claims in the theorem. The first lemma establishes the FSLLN.

\begin{lemma}[FSLLN]\label{lem:arrival-fluid}
The multivariate traffic process $\mathbf A_n = (A_1, \ldots, A_J) :=
\sum_{m=1}^n \mathbf a_m$ satisfies
\[
n^{-1} \mathbf A_n \to \mathbf F~\text{in}~(\sC^J,U)~\text{a.s.}
\]
as $n \to \infty$, where $\mathbf F = (F_1,\ldots, F_J)$ and $F_j(t) =
\bbE[\mathbf 1_{\{T_{j} \leq t\}}]$ for all $t \in [0,T]$.
\end{lemma}
\Proof
  First, for each $j \in \sE$, the classical Glivenko-Cantelli theorem
  implies that
  \begin{equation}
    \label{eq:indi-fluid}
    n^{-1} A_{j} \to F_j ~\text{in}~(\sC,U)~\text{a.s.} 
  \end{equation}
  as $n \to \infty$. By the multivariate strong law of large numbers
  it follows that for a fixed $t \in [0,T]$
\(
    \mathbf A_n(t) \to \mathbf F(t)~\text{a.s.}
\)
  as $n \to \infty$. The functional limit follows as a consequence
  of~\eqref{eq:indi-fluid}.
\EndProof

This proves part $(i)$ of Theorem~\ref{thm:arrival-limits}. The next
lemma establishes part $(ii)$.

\begin{lemma}
  The multivariate traffic process $\mathbf A_n$ satisfies a
  functional central limit theorem where
\[
\sqrt n \left(n^{-1} \mathbf A_n - \mathbf F \right) \Rightarrow W^0
\circ \mathbf F~\text{in}~(\sC^J,U),
\]
where $\mathbf W^0 \circ \mathbf F$ is a $J$-dimensional Brownian
bridge process as defined in Definition~\ref{def:bb-multi}, with
covariance function $(R(t),~t\geq 0) = ([F_{i,j}(t) - F_i(t)F_j(t)],~t
\geq 0)$.
\end{lemma}
\Proof
Once again, Donsker's theorem for empirical processes implies that 
\begin{equation}\label{eq:indi-diff}
  \hat A_j  := \sqrt n \left( n^{-1} A_j - F_j \right) \Rightarrow W_j^0 \circ F_j ~\text{in}~ (\sC,U)
\end{equation}
as $n \to \infty$ for every $j \in \sK$. This implies that the
marginal arrival processes are tight. \cite[Theorem 11.6.7]{Wh01}
implies that the multivariate process $\mathbf A_n$ is also
tight.The multivariate central limit theorem
\cite[Theorem 4.3.4]{Wh01} implies that the scaled process $\hat
{\mathbf A}_n(t) = (\hat A_1(t), \ldots, \hat A_J(t))$ (for fixed $t
\in [0,T]$) satisfies 
\[
\hat{\mathbf A}_n (t) = \sqrt n \left(\frac{\mathbf A_n(t)}{n} -
  \mathbf F(t)\right) \Rightarrow \mathcal N(0,R(t)),
\]
where $\mathcal N(0,R(t))$ is a mean zero $J$-dimensional Gaussian random vector
with covariance matrix $R(t) = [F_{i,j}(t) - F_i(t)F_j(t)]$. The Cram\'er-Wold device together
with this result implies that the finite-dimensional distributions of
$\mathbf A_n$ converge weakly to a tuple of Gaussian random
vectors. The tightness of the processes $\{\mathbf A_n\}$, the
continuity of the limit process  and Prokhorov's
theorem implies that $\hat{\mathbf{A}}_n$ converges weakly to
the multivariate Gaussian stochastic process $\mathbf W^0 \circ \mathbf F$ with mean zero and
covariance function $(R(t), t\geq 0)$ in $(\sC^J, U)$.
\EndProof

% \subsection{Proof of Theorem~\ref{thm:service-limits}}
% As in the previous section, we divide the proof into two lemmas, one
% focusing on the fluid limit and the second on the diffusion
% limit. Recall that we assume the existence of a sequence of
% non-decreasing continuous functions $\{\mathbf M_n,~n\geq 1\}$ where
% $\mathbf M_n \in \sD^K^K$ and that for each $k$ there exists $\mu_{n,k}
% \in \sD^K$ non-negative such that $M_n(t) = \int_0^t \mu_{n,k}(s)
% ds$. We also assume that there exists $\mathbf M \in \sD^K$ such that
% $\mathbf M_n \to \mathbf M$ in $(\sC^K,U)$ as $n \to infty$.

% \begin{lemma}
%   Let $\{\mathbf S_n := (S_{n,1},\ldots, S_{n,k}),~n\geq 1 \}$ be a
%   sequence of counting processes, with independent component processes
%   $S_{n,k}$ for $k \in \sK$. Suppose that $\|\bbE[\mathbf S_n] - M_n\|
%   \to 0$ as $n \to \infty$. Then, 
% \[
% \frac{\mathbf S_n}{n} \to \mathbf M ~\text{in}~(\sC^K,U)~\text{a.s. as}
% $n \to \infty$.
% \]
% \end{lemma}
\bibliographystyle{plainnat} 
\bibliography{refs-queueing}

\end{document}